\documentclass[10pt]{article}
\newcommand{\documentdate}{29 IX 2026}
\usepackage{a4wide,latexsym,graphicx,amsmath,amssymb,varioref}
\usepackage{graphics,xcolor,dsfont}
\usepackage{cancel}
\usepackage{ulem}

\newcommand{\al}[1]{{\footnotesize{\sf #1}}}
\title{{\sf ADSWITCH}: A Simple First-Order Algorithm \\for Full-Rank Equality-Constrained Optimization}

\author{
S. Gratton\thanks{Universit\'{e} de Toulouse, INP, IRIT, Toulouse, France. Email:
     serge.gratton@enseeiht.fr. Work partially supported by 3IA Artificial and
     Natural Intelligence Toulouse Institute (ANITI), French "Investing for the Future
     - PIA3" program under the Grant agreement ANR-19-PI3A-0004"}
~and Ph. L. Toint\thanks{Universit\'{e} de Toulouse, INP, IRIT,
  Toulouse, France, and NAXYS, University of Namur, Namur, Belgium. Email: philippe.toint@unamur.be}
}

\DeclareMathAlphabet{\pazocal}{OMS}{zplm}{m}{n}

\newcommand{\calO}{{\pazocal{O}}}

\newcommand{\beqn}[1]{\begin{equation}\label{#1}}
\newcommand{\eeqn}{\end{equation}}
\newcommand{\req}[1]{(\ref{#1})}
\newcommand{\ms}{\;\;\;\;}
\newcommand{\tim}[1]{\;\; \mbox{#1} \;\;}

\newtheorem{theorem}{Theorem}[section]
\newtheorem{lemma}[theorem]{Lemma}

\newcommand{\numsection}[1]{\section{#1}\setcounter{equation}{0}}
\newtheorem{corollary}{Corollary}
\newcommand{\appnumsection}[1]{\section*{#1}\setcounter{equation}{0}
  \renewcommand{\theequation}{A.\arabic{equation}}
  \renewcommand{\thetheorem}{A.\arabic{theorem}}
  \renewcommand{\thetable}{A.\arabic{table}}
  \renewcommand{\thefigure}{A.\arabic{figure}}
  \renewcommand{\thesection}{A} }
\renewcommand{\theequation}{\arabic{section}.\arabic{equation}}
\renewcommand{\thefootnote}{(\arabic{footnote})}
\newcounter{algo}[section]
\renewcommand{\thealgo}{\thesection.\arabic{algo}}
\newcommand{\llem}[2]{\vspace{\baselineskip} 
\noindent\framebox[\textwidth]{\parbox{0.95\textwidth}{
\begin{lemma} \label{#1} \rm #2 \end{lemma} } } \vspace{\baselineskip} }

\newcommand{\algo}[3]{\refstepcounter{algo}
\begin{center}\begin{figure}[htbp]
\framebox[\textwidth]{
\parbox{0.95\textwidth} {\vspace{\topsep}
{\bf Algorithm \thealgo : #2}\label{#1}\\
\vspace*{-\topsep} \mbox{ }\\
{#3} \vspace{\topsep} }}
\end{figure}\end{center}}
\newcommand{\bpr}{{\bf Proof.} \hspace{1.5mm}}
\newcommand{\epr}{\hfill $\Box$ \vspace*{1em}}
\newcommand{\proof}[1]{
\begin{list}{}{
\setlength{\topsep}{0.0pt}
\setlength{\partopsep}{0.0pt}
\setlength{\leftmargin}{0.025\textwidth}
\setlength{\rightmargin}{0.5\leftmargin}
\setlength{\labelwidth}{0.5\leftmargin}
\setlength{\labelsep}{0.25\leftmargin}}
\item \bpr #1 \epr \noindent
\end{list}}
\newcommand{\lthm}[2]{\vspace{\baselineskip} 
\noindent\framebox[\textwidth]{\parbox{0.95\textwidth}{
\begin{theorem} \label{#1} \rm #2 \end{theorem} } } \vspace{\baselineskip} }

\newcommand{\iibe}[2]{\{ #1, \ldots, #2 \}}
\renewcommand{\Re}{\hbox{I\hskip -2pt R}}
\newcommand{\smallRe}{\hbox{\footnotesize I\hskip -2pt R}}

\newcommand{\sfrac}[2]{{\scriptstyle \frac{#1}{#2}}}
\newcommand{\kap}[1]{\kappa_{\mbox{\tiny #1}}}
\newcommand{\eqdef}{\stackrel{\rm def}{=}}
\newcommand{\tal}[1]{{\normalsize {\sf #1}}}
\newcommand{\half}{\sfrac{1}{2}}
\newcommand{\flow}{f_{\rm low}}

\newcommand{\comment}[1]{}
\newcommand{\private}[1]{}
\date{\documentdate}

\begin{document}

\renewcommand{\thefootnote}{\fnsymbol{footnote}}
\maketitle
\renewcommand{\thefootnote}{\arabic{footnote}}

\abstract{
A very simple first-order algorithm is proposed for solving
nonlinear optimization problems with  nonlinear equality
constraints. The algorithm adaptively selects steps in the plane
tangent to the constraints or steps that reduce infeasibility, without
using a merit function or filter.  The tangent steps are based on the
AdaGrad method for unconstrained minimization. The objective function
is never evaluated by the algorithm, making it suitable for noisy
problems.  Its worst-case evaluation complexity is analyzed, yielding
a global convergence rate of $\calO(1/\sqrt{k})$, where $k$ is the
iteration number. This order matches the best-known rate of
first-order methods for unconstrained problems. Numerical experiments
are presented suggesting that the performance of the algorithm is
comparable to that of first-order methods for unconstrained problems,
and that its reliability is remarkably stable in the presence of noise
in the gradient of the objective function.
}


\numsection{Introduction}

Solving nonlinear problems with nonlinear equality constraints is a
central part of more general constrained optimization, and the purpose
of this paper is to contribute to the subject by proposing a very
simple algorithm for the case where the constraints' Jacobian has
full rank.  We introduce a method that draws its inspiration from both
the ``trust-funnel'' approaches \cite{GoulToin10,CurtRobiSama18b}
and the objective-function-free (OFFO) first-order algorithms that
have become very popular because of their successful use in deep
learning applications. Trust-funnel methods themselves exploit the
older idea of decomposing the step computed within an optimization
algorithm for constrained problems into ``tangential'' and
``normal'' components. Advocated in particular by Byrd and Omojokun
\cite{Omoj89} for sequential quadratic programming, this last idea is
to consider that a useful minimization step should, on the one hand,
improve the value of the objective function without increasing the
constraints' violation too much (this is the tangential step, because
it can be interpreted as a step (nearly) tangent to the constraint's
surface) and, on the other hand, improve feasibility (the normal step,
mostly orthogonal to the constraints' surface). This technique has
been successfully used by many algorithms, and in particular by
trust-funnel ones.  In this case, the two types of step are treated
separately and a choice between normal and tangential steps is made at
every iteration by comparing constraint violation (primal feasibility)
with some dual optimality measure. Our proposed algorithm uses a
similar mechanism, although considerably simplified.  Our second
source of inspiration is the thriving field of first-order methods for
unconstrained problems, and, more specifically, the well known and
well understood AdaGrad algorithm
\cite{DuchHazaSing11,WardWuBott19,DefoBottBachUsun22,GratJeraToin22b}
or unconstrained optimization. Such methods are attractive because they avoid computing
the value of the objective function entirely (they only rely on
gradients). This feature makes them quite robust in the presence of
noise in the gradient of the objective function, as is commonly the case when its evaluation involves
subsampling.  Again, our method exploits this
OFFO\footnote{Objective-Function-Free Optimization.} strategy to
achieve robustness.  In effect, it consists of a simplified
trust-funnel method in which the tangential step is computed using the
AdaGrad algorithm.

Other OFFO methods for constrained optimization have been proposed.
The stochastic algorithm proposed in \cite{BeraCurtRobiZhou21} uses a combined
tangential and normal step with a quadratic model of a merit function. A
(nontrivial) stochastic analysis of its convergence (albeit not its worst-case
complexity) is provided. It also allows the use of second-order information, at the
cost of a somewhat complicated algorithm. \cite{CurtRobiZhou24} extends
this proposal to handle inequality constraints and provides an updated
convergence analysis.
The algorithm of \cite{FangNaMahoKola24} is more similar to ours, but differs
in significant aspects. Although it uses tangential and normal steps,
it does so at every iteration. The step in the Jacobian's nullspace is
computed by a standard trust-region technique and a special rule is
provided to update its radius. It also relies on an explicit penalty
parameter, which is updated in the course of the algorithm. Convergence
is analyzed in a stochastic setting.
The method of \cite{CoelFernCostFerr24} is inspired by filter methods
(another class of methods for nonlinear optimization
\cite{FletLeyf02,FletLeyfToin02,FletGoulLeyfToinWaec02}) and also alternates
between steps that attempt to reduce constraint infeasibilities and
``optimality'' steps improving the objective function value,
albeit ignoring  constraints in the second case. Good numerical
results are reported on deep learning problems, but no convergence
theory or complexity analysis is provided\footnote{Since the
optimality steps ignore the constraints geometry, it seems possible that
they might significantly decrease feasibility and that cycling
could occur.}.

In view of these comments, we summarize our contributions as follows.
\begin{enumerate}
\item In order to solve problem \req{problem}, we propose an
  ``adaptive switching'' algorithm, dubbed \al{ADSWITCH}, which can be 
  viewed as a simplified trust-funnel method. \al{ADSWITCH} uses the
  purely OFFO first-order AdaGrad algorithm in the plane tangent to
  the constraints. 
\item We provide a complete worst-case complexity analysis for the
  case where the objective's gradient and the constraints are 
  deterministic, under assumptions similar to those used in
  \cite{BeraCurtRobiZhou21,CurtRobiZhou24}.
\item We illustrate the practical behaviour of the method on a subset
  of the {\sf CUTEst} problems \cite{GoulOrbaToin15b} as provided by
  the S2MPJ environment \cite{GratToin24}.
\end{enumerate}

The paper is organized as follows. The \al{ADSWITCH} algorithm is
introduced in Section~\ref{sec:algo} together with the necessary
concepts and notation.  Section~\ref{sec:det-comp} contains the
relevant complexity analysis. Numerical results are proposed in
Section~\ref{sec:numerics}, while some conclusions and perspectives are
discussed in Section~\ref{sec:concl}.

\noindent
{\bf Notation:} In what follows, $\|\cdot\|$ denotes the Euclidean
norm of a vector in the relevant space, and $\sigma_{\min}(M)$ is the smallest
singular value of the matrix $M$. 

\numsection{An adaptative switching algorithm}\label{sec:algo}

We propose solving the smooth equality constrained problem
\beqn{problem}
\min_{x\in\smallRe^n} f(x)\qquad \text{such that}\qquad c(x)=0,
\eeqn 
using the purely first-order \al{ADSWITCH} method, which we now 
describe after establishing some notation. In \req{problem}, $f$ is a
smooth function from $\Re^n$ into $\Re$ and $c$ is a smooth function
from $\Re^n$ into $\Re^m$ ($m \le n$). For a given vector $x$, we
assume that we can compute  the
gradient $g(x) = \nabla_x f(x)$, as well as the (exact) Jacobian of
the constraints at $x$,
\[
J(x)=\nabla c(x)\in\Re^{m\times n}.
\]
Assuming $J(x)$ has full rank, we may also compute the orthogonal projection onto its nullspace by
\beqn{proj-def}
P(x)\eqdef I - J(x)^T\big(J(x)J(x)^T\big)^{-1}J(x).
\eeqn
The projected gradient is 
\beqn{gT-def}
g_T(x) \eqdef P(x)\,g(x).
\eeqn
We may also define the (exact) least-squares Lagrange multiplier
$\widehat\lambda(x)$ by
\beqn{eq:ls-mult}
\big(J(x)J(x)^T\big)\,\widehat\lambda(x)\ =\ -\,J(x)\,\nabla f(x),
\eeqn
the standard Lagrangian
\beqn{Lag-def}
L(x,\lambda) = f(x) + \lambda^Tc(x)
\eeqn
and, for $\rho\ge 0$ fixed, the associated ``augmented-Lagrangian-like'' Lyapunov function
\beqn{lyap-def}
\psi_\rho(x,\lambda) \eqdef L(x,\lambda)+\rho\|c(x)\|,
\eeqn
which we often abbreviate for $\lambda = \widehat\lambda(x)$ as
\beqn{psi-def}
\psi(x)\eqdef \psi_\rho\big(x,\widehat\lambda(x)\big).
\eeqn
The definitions imply the important properties
\beqn{eq:ids}
g_T(x)=\nabla f(x)+J(x)^T \widehat\lambda(x)
\tim{ and }
J(x)\,g_T(x)=0
\eeqn
and also, for $\lambda$ fixed to $\widehat\lambda(x)$,
\beqn{eq:gradL}
\nabla_x L\big(x,\lambda\big)= g_T(x).
\eeqn
By construction, we also have
\beqn{JTc-gT}
c(x)^TJ(x)g_T(x)=0.
\eeqn
Our objective is now to compute a first-order critical point $x_*$ for
problem \req{problem}, that is a point $x_*$ such that
\[
g_T(x_*) = P(x_*)g(x_*) = 0
\tim{ and }
c(x_*) = 0.
\]
(Observe that the objective function does not appear in these
conditions.)

We are now in a position to specify \al{ADSWITCH} in detail.
As announced in the introduction, the algorithm
uses two types of step: ``tangential steps'' attempt to reduce the
objective-function value (without evaluating it!) in the nullspace of
the constraints, while ``normal steps'' aim at reducing infeasibility.
Which type of step is used at a given iteration depends on the
respective values of the tangential step and constraint violation.
 
\algo{adagec}{\tal{ADSWITCH}}{
  \begin{description}
  \item[Step 0: Initialization.]  A starting point $x_0$ is
    given, together with constants  $\theta>1$,
    $\beta,\eta,\varsigma\in (0,1]$. Set $\Gamma_0= 0$ and $k=0$.
  \item[Step 1: Evaluations.]
    Evaluate $c_k=c(x_k)$, $J_k = J(x_k)$ and $g_k= \nabla_xf(x_k)$.\\
    Set $g_{T,k} = P(x_k)g_k$,
    \beqn{Gamma-def}
    \Gamma_k^+ = \Gamma_k + \|g_{T,k}\|^2,
    \eeqn
    and
    \beqn{alphaT-def}
    \alpha_{T,k} = \frac{\eta}{\sqrt{\Gamma_k^+ + \varsigma}}.
    \eeqn
  \item[Step 2: Tangential step: ]
    If
    \beqn{switch}
    \|c_k\| \leq \beta \alpha_{T,k}\|g_{T,k}\|
    \eeqn
    then set
    \beqn{xplusT-def}
    x_{k+1}=x_k- \alpha_{T,k}g_{T,k}
    \eeqn
    and
    \beqn{Gamma-upd}
    \Gamma_{k+1} = \Gamma_k^+,
    \eeqn
    and go to Step~4.
  \item[Step 3: Normal step:]
    Otherwise, find a step  $s_{N,k}$ in the range space of $J_k^T$
    such that
    \beqn{TR-bound}
    \|s_{N,k} \| \leq \theta \|c_k\|
    \eeqn
    and there exists $\kap{nrm}>0$ independent of $k$ such that
    \beqn{normal-descent}
    \half \|c(x_k+s_{N,k})\|^2 \le  \half\|c_k\|^2 - \kap{nrm} \|J_k^Tc_k\|^2,
    \eeqn
    and set
    \beqn{xplusN-def}
    x_{k+1} = x_k + s_{N,k} \tim{ and } \Gamma_{k+1}=\Gamma_k.
    \eeqn
  \item[Step 4: Loop.]
    Increment $k$ by one and go to Step~1.
\end{description}
}

\begin{itemize}
\item The tangential step is the standard AdaGrad-norm
  \cite{WardWuBott19} step in the  nullspace of the Jacobian.
\item The statement of the procedure to find the normal step in Step~3
  may seem abstract because $\kap{nrm}$ may not be known, but,
  fortunately, there exist several standard techniques to achieve
  \req{TR-bound}--\req{normal-descent}.
  
  A first technique is to use a standard Armijo backtracking
  linesearch along the steepest-descent direction $-J_k^Tc_k$ (see
  Lemma~\ref{csq-decrease}).  It is well known that
  \req{normal-descent} holds in this case with $\kap{nrm}= 1/2L_{JTc}$
  (see \cite[Lemma~2.2.1]{CartGoulToin22}, for instance).  More
  generally, any ``gradient related'' direction $d_k$, that is any
  $d_k$ satisfying \req{TR-bound} and $d_k^TJ_k^Tc_k\leq
  -\kap{grl}\|J_k^Tc_k\|^2$ for some $\kap{grl}>0$, can be used to
  initiate a backtracking linesearch on $\half\|c(x)\|^2$.
  
  A second possibility follows along this line and is to select
  a positive definite matrix $B_k$ with bounded condition number to
  define $d_k=-B_k^{-1}J_k^Tc_k$ and then to perform a backtracking
  linesearch along $d_k$ (see \cite[Section~10.3.1]{ConnGoulToin00}
  for a proof that $d_k$ is a gradient-related direction), in which
  case a similar result holds.  For instance, using
  $B_k = J_k^TJ_k + \delta I_n$ for some $\delta \ge 0$ yields the
  regularized Gauss-Newton step
  \beqn{reg-Newton-normal}
  s_{Nk} = - \gamma_k (J_k^TJ_k + \delta I_n)^{-1}J_k^Tc_k
  = -\gamma_k J_k^T (J_k J_k^T + \delta I_m)^{-1}c_k,
  \eeqn
  where $\gamma_k>0$ is chosen by backtracking as large as possible to ensure
  \req{TR-bound} and \req{normal-descent}. (Note that, for $\delta>0$,
  $J_k^TJ_k + \delta I_n$ has bounded condition number when $J_k$ is bounded.)

  Suitable variants of trust-region steps are also possible
  \cite[Section~10.3.2]{ConnGoulToin00}, themselves involving several
  subvariants such as the truncated conjugate gradient
  \cite{Stei83a,Toin81b}, CGLS \cite{HestStie52} or LSQR
  \cite{PaigSaun82a}.  Because iterates generated by these methods are
  linear combinations of the gradients of $\half\|c_k+J_k s\|^2$, t the final
  is guaranteed to belong to the range of $J_k^T$.
  
\item We have mentioned an explicit formula for the projection
  $P(x)$ in \req{proj-def}, but it should be computed using a rank-revealing QR
  factorization of $J_k^T$:
  \[
  Q_k^TJ_k^TP_k
  = \left(\begin{array}{c}Y_k^T\\Z_k^T\end{array}\right)J_k^TP_k
  =\left(\begin{array}{c}R_k\\0\end{array}\right),
    \]
  where $P_k$ is a permutation matrix and the columns of $Z_k$
  provide an orthonormal basis for the nullspace of $J_k$. 
  In this case,
  \[
  g_{T,k} = Z_kZ_k^T g_k
  \]
  and the Newton step \req{reg-Newton-normal} is given by
  \beqn{facto-Newton-normal}
  s_{N,k} = \left\{\begin{array}{ll}
  -\gamma_kY_kR_k^{-T}P_kc & \tim{if }\delta=0,\\
  -\gamma_k Y_kR_k ( R_k^T R_k + \delta I_m)^{-1}P_k c_k
  & \tim{otherwise.}
  \end{array}\right.
  \eeqn
\item Observe that \al{ADSWITCH} does not involve
  any explicit merit function or filter to control overall progress,
  but relies on the simple switching condition \req{switch} (hence its
  name). As we see below, an implicit merit function (a modified
  augmented Lagrangian with suitably chosen parameters) is nevertheless
  crucial for our theoretical argument.
\item \al{ADSWITCH} differs in several aspects from
  trust-funnel methods of \cite{GoulToin10,CurtRobiSama18b}.
  \begin{enumerate}
    \item It avoids evaluating the objective function in order
      to improve its reliability in the presence of noise in the
      objective function's gradient.
    \item It also avoids estimating Lagrange multipliers, but
      at the cost of computing the exact projection $P(\cdot)$ defined
      in \req{proj-def}.
    \item The trust-region approach, which is explicit in
  trust-funnel methods, is only implicit in \al{ADSWITCH}. Indeed it
  has been argued in \cite{GratJeraToin22a} that the (unconstrained)
  AdaGrad algorithm can be viewed as a trust-region method with a
  specific radius management. Moreover, \req{TR-bound} can
  be viewed as a trust-region constraint (such a condition is
  also imposed in \cite{GoulToin10}), but without the need to
  explicitly update its radius. \al{ADSWITCH}, in its
  present incarnation, is more limited than trust-funnel methods in
  that it does not handle approximate projections in the tangent
  plane.
  \end{enumerate}
\end{itemize}

\noindent
In what follows, we denote by $\{k_\tau\}\subseteq \{k\}$ the index
  subsequence of ``tangential iterations'' where
  \req{switch} holds. We also denote by $\{k_\nu\} =
    \{k\}\setminus\{k_\tau\}$ the index subsequence of ``normal iterations''.
By convention, sums over empty sets of indices are zero.
    
\numsection{Complexity analysis}\label{sec:det-comp}

\subsection{Assumptions}

Our analysis uses the following assumptions.

\begin{description}
\item[A0:] $f$ and $c$ are continuously differentiable on $\Re^n$.
\item[A1:] For all $x\in\Re^n$, $f(x)  \geq \flow$.
\item[A2:] For all $x\in\Re^n$, $\|g(x)\| \leq \kappa_g$ where $\kappa_g\ge \eta\beta$.
\item[A3:] For all $x\in\Re^n$, $\|J(x)\| \leq \kappa_J$.
\item[A4:] For all $x\in\Re^n$, $\sigma_{\min}(J(x)) \geq \sigma_0 \in (0,1]$.
\item[A5:] The gradient $g(x)$ is globally Lipschitz continuous (with constant $L_g$).
\item[A6:] The Jacobian $J(x)$ is globally Lipschitz continuous (with constant $L_J$).
\end{description}
Observe that A1, A2 and A3 are automatically satisfied if
the iterates $x_k$ remain in a bounded domain of $\Re^n$. The
following results also follows from these assumptions.

\llem{c-bounded}{Suppose that A0, A2, A3 and A6 hold. Then there exists a
  constant $\kappa_c>0$ such that $\|c_k\|\le \kappa_c$ for all
  $k\ge 0$.
}
\proof{
We have from the definition of $g_{T,k}$, the contractive nature of
the projection and A2 that $\|g_{T,k}\|\le \|g_k\| \le \kappa_g$.
Thus, \req{alphaT-def} ensures that, at a tangential iteration, the step
$s_{T,k}$ satisfies
\[
\|s_{T,k}\|
= \alpha_{T,k}\|g_{T,k}\|
\leq \frac{\eta\kappa_g}{\varsigma}
\eqdef \kap{sw},
\]
The switching condition \req{switch} gives
$\|c_k\|\leq\beta\alpha_{T,k}\|g_{T,k}\|\leq\beta\kap{sw}$.
Using A3 and A6 to bound the remainder in the first-order expansion
of $c$, we therefore obtain
\[
\|c_{k+1}\|
\leq \|c_k\|+\|J_k\|\,\|s_{T,k}\|+\frac{L_J}{2}\|s_{T,k}\|^2
\leq \beta\kap{sw}+\kappa_J\kap{sw}+\frac{L_J}{2}\kap{sw}^2.
\]
At a normal iteration, \req{normal-descent} gives
$\|c_{k+1}\|\leq\|c_k\|$. The desired bound follows by induction
with $\kappa_c =
\max[\|c(x_0)\|,(\beta+\kappa_J)\kap{sw}+\frac{L_J}{2}\kap{sw}^2]$
\epr
}

\noindent
Also note that
\begin{itemize}
\item A4 ensures that $P(x_k)$,  $g_{T,k}$ and $\widehat\lambda(x)$ are
well defined;

\item A3, A5 and A6 ensure that
$\widehat\lambda(x)$ is Lipschitz continuous (with constant $L_\lambda$);

\item A3 ensures that $c$ is Lipschitz continuous (with constant $L_c$);
  
\item A3, A4, A5 and A6  ensure that
  $\nabla_x L(x,\widehat\lambda)$ is Lipschitz
  continuous (with constant $L_L$) as a function of $x$ for fixed $\widehat\lambda$.
\end{itemize}
These assumptions are very similar to those of
\cite{BeraCurtRobiZhou21,CurtRobiZhou24}. In particular, A4 is part of
Assumption~1 in this last reference.

\subsection{Normal steps}

Our analysis hinges on the fact that first-order descent can be shown
on the Lyapunov $\psi(x)$, for both tangential and normal
steps, despite the fact that neither $\widehat\lambda(x_k)$ nor $\rho$
(which we still need to define) appears in the algorithm.
We start by considering normal steps for a sufficiently large penalty
parameter $\rho$.
We start with a fairly simple observation.

\llem{csq-decrease}{Suppose A4 holds and a normal step is used at
  iteration $k_\nu$.  Then
  \beqn{c-descent}
  \|c_{k_\nu+1}\| \le \left(1-\kap{nrm}
  \sigma_0^2 \right)\|c_{k_\nu}\|.
  \eeqn
}

\proof{
  We have from \req{normal-descent} that
  $\|c_{k_\nu+1}\|<\|c_{k_\nu}\|$.  Then
  \[
  2\|c_{k_\nu}\|(\|c_{k_\nu+1}\|-\|c_{k_\nu}\|)
  \le (\|c_{k_\nu+1}\|+\|c_{k_\nu}\|)(\|c_{k_\nu+1}\|-\|c_{k_\nu}\|)
  = \|c_{k_\nu+1}\|^2-\|c_{k_\nu}\|^2,
  \]
  and therefore, using \req{normal-descent} again, 
  \[
  \|c_{k_\nu+1}\|-\|c_{k_\nu}\|
  \le -\frac{\kap{nrm}\|J_{k_\nu}^Tc_{k_\nu}\|^2}{\|c_{k_\nu}\|}
  \le -\frac{\kap{nrm}\sigma_0^2\|c_{k_\nu}\|^2}{\|c_{k_\nu}\|}
  = -\kap{nrm}\sigma_0^2\|c_{k_\nu}\|.
  \]
} \epr

\llem{lem:normal-descent}{
Suppose A0--A6 hold and a normal step is used at
iteration $k_\nu$. Define
\beqn{rho-def}
\rho = \frac{1}{\kappa_{\rm nrm}\,\sigma_0^2}\Big[\;\kappa_c\,\theta\,\,L_\lambda\;+\;\frac{L_L}{2}\,\kappa_c\,\theta^2\;+\;\eta\;\Big].
\eeqn
Then we have
\beqn{effectN}
\psi(x_{k_\nu+1}) - \psi(x_{k_\nu}) \leq -\eta \|c_{k_\nu}\|.
\eeqn
}

\proof{
We have
\beqn{DaDb}
\begin{aligned}
\psi(x_{k_\tau+1}) - \psi(x_{k_\tau})
&= \underbrace{\psi_\rho(x_{k_\tau+1},\widehat\lambda(x_{k_\tau}))-\psi_\rho(x_{k_\tau},\widehat\lambda(x_{k_\tau}))}_{\Delta_x}\\
&\ms\ms\ms +\underbrace{\psi_\rho(x_{k_\tau+1},\widehat\lambda(x_{k_\tau+1}))-\psi_\rho(x_{k_\tau+1},\widehat\lambda(x_{k_\tau}))}_{\Delta_\lambda}.
\end{aligned}
\eeqn
We consider $\Delta_x$ and $\Delta_\lambda$ separetely.
The first part of
\req{eq:ids} and the fact that $s_{N,k}$ belongs to the orthogonal
complement of
the Jacobian's nullspace, $\nabla_xL(x_{k_\nu},\widehat{\lambda}_{k_\nu})=g_{T,k_\nu}$ and
$g_{T,k_\nu}^Ts_{N,k_\nu}=0$. Thus, using \req{psi-def},
\req{c-descent}, \req{TR-bound}, the Lipschitz continuity of
$\nabla_xL(x,\widehat\lambda)$ for fixed $\widehat\lambda$ ($\rho$ is
fixed in \req{rho-def}) and the definition of $\Delta_x$ and 
$\Delta_\lambda$ in \req{DaDb}, we obtain
\beqn{eq:N1}
\begin{aligned}
\Delta_x
& = \psi_\rho(x_{k_\nu+1},\widehat\lambda_{k_\nu})-\psi_\rho(x_{k_\nu},\widehat\lambda_{k_\nu})\\
& = L(x_{k_\nu+1},\widehat{\lambda}_{k_\nu})-L(x_{k_\nu},\widehat{\lambda}_{k_\nu})
+\rho\big(\|c_{k_\nu+1}\| - \|c_{k_\nu}\|\big)\\
& \leq \nabla_xL(x_{k_\nu},\widehat\lambda_{k_\nu})^Ts_{N,k_\nu} + \frac{L_L}{2}\|s_{N,k_\nu}\|^2
- \rho \kap{nrm}\sigma_0^2\|c_{k_\nu}\|\\
& \leq - \rho \kap{nrm}\sigma_0^2\|c_{k_\nu}\| + \frac{\theta^2L_L}{2}\|c_{k_\nu}\|^2.
\end{aligned}
\eeqn
Using the Lipschitz continuity of $\widehat\lambda$, \req{TR-bound} and \req{normal-descent},
we then deduce that
\beqn{eq:N2}
\begin{aligned}
\Delta_\lambda
& = \psi_\rho(x_{k_\nu+1},\widehat\lambda_{k_\nu+1})-\psi_\rho(x_{k_\nu+1},\widehat\lambda_{k_\nu})\\
&\le \|c_{k_\nu+1}\|\,\|\widehat\lambda_{k_\nu+1}- \widehat\lambda_{k_\nu}\|\\
&\le L_\lambda\,\|c_{k_\nu+1}\|\,\|s_{N,k_\nu}\|\\
&\le L_\lambda\,\theta \|c_{k_\nu}\|^2, \\
\end{aligned}
\eeqn
and thus, summing \req{eq:N1} and \req{eq:N2} and using Lemma~\ref{c-bounded}, that
\[
\psi(x_{k_\nu+1}) - \psi(x_{k_\nu})
\leq - \rho\kap{nrm}\sigma_0^2\|c_{k_\nu}\|
+L_\lambda\,\theta\kappa_c\,\|c_{k_\nu}\|
+\frac{L_L}{2}\theta^2\kappa_c\|c_{k_\nu}\|.
\]
The bound \req{effectN} then follows from \req{rho-def}.
} \epr

\subsection{Tangential steps}

Having considered normal steps, we now show that a decrease in the
Lyapunov function $\psi$ may also be proved if a tangential step is taken.

\llem{tangent-decrease}{
Suppose A0--A6 hold. Then
\beqn{effectT}
\psi(x_{k_\tau+1}) - \psi(x_{k_\tau})
\le -\alpha_{T,k_\tau}\|g_{T,k_\tau}\|^2
+\kap{tan}\,\alpha_{T,k_\tau}^2\|g_{T,k_\tau}\|^2,
\eeqn
where
\[
\kap{tan} = \frac{1}{2}\Big(L_L+\rho L_J\Big)  + \beta L_\lambda  + L_\lambda L_c.
\]
}

\proof{We again consider separately the term $\Delta_x$ and $\Delta_\lambda$  in the
decomposition \req{DaDb}. The Lipschitz
continuity of  $\nabla_xL(x,\widehat\lambda)$ for fixed
$\widehat\lambda$ and \req{eq:gradL} give
\[
\Delta_x = -\nabla_x L(x_{k_\tau},\widehat\lambda(x_{k_\tau}))^T(\alpha_{T,k_\tau}g_{T,k_\tau})
  + r_0 + \rho(\|c(x_{k_\tau+1})\|-\|c(x_{k_\tau})\|),
\]
with
\[
|r_0| \leq \frac{L_L}{2}\alpha_{T,k_\tau}^2\|g_{T,k_\tau}\|^2.
\]
Now, using A6, 
\[
\|c(x_{k_\tau+1})\| =
\|c(x_{k_\tau})-\alpha_{T,k_\tau}J_kg_{T,k_\tau} + r_1\|
\le \|c(x_{k_\tau})\| + \|r_1\|
\]
with
\[
\|r_1\| \leq \frac{L_J}{2}\alpha_{T,k_\tau}^2\|g_{T,k_\tau}\|^2,
\]
giving
\beqn{Encplus}
\|c(x_{k_\tau+1})\|-\|c(x_{k_\tau})\|
\le \frac{L_J}{2}\alpha_{T,k_\tau}^2\|g_{T,k_\tau}\|^2.
\eeqn
Successively using \req{eq:gradL} and \req{JTc-gT}, we obtain
\beqn{desc-det}
\begin{aligned}
\Delta_x 
&\le -\alpha_{T,k_\tau}g_{T,k_\tau}^Tg_{T,k_\tau} + r_0
+ \frac{\rho L_J}{2}\,\alpha_{T,k_\tau}^2\|g_{T,k_\tau}\|^2\\
&\le -\alpha_{T,k_\tau}g_{T,k_\tau}^Tg_{T,k_\tau} +\frac{1}{2}\left(L_L+\rho L_J\right)\,\alpha_{T,k_\tau}^2\|g_{T,k_\tau}\|^2 \\
&\le -\alpha_{T,k_\tau}\|g_{T,k_\tau}\|^2
+\frac{1}{2}\left(L_L+\rho L_J\right)\,\alpha_{T,k_\tau}^2\|g_{T,k_\tau}\|^2.
\end{aligned}
\eeqn
Now, we may use the Lipschitz continuity of $\widehat\lambda$ and $c$ to deduce that
\beqn{Dlprod}
\begin{aligned}
\Delta_\lambda
& =c_{k_\tau+1}^T\big(\widehat\lambda_{k_\tau+1}-\widehat\lambda_{k_\tau}\big)\\
&= (c_{k_\tau+1}-c_{k_\tau})^T\big(\widehat\lambda_{k_\tau+1}-\widehat\lambda_{k_\tau}\big)
 + c_{k_\tau}^T\big(\widehat\lambda_{k_\tau+1}-\widehat\lambda_{k_\tau}\big)\\
&\le \|c_{k_\tau}\|\,\|\widehat\lambda_{k_\tau+1}-\widehat\lambda_{k_\tau}\|
 + \|c_{k_\tau+1}-c_{k_\tau}\|\,\|\widehat\lambda_{k_\tau+1}-\widehat\lambda_{k_\tau}\|\\
&\le L_\lambda \|c_{k_\tau}\|\,\alpha_{T,k_\tau}\|g_{T,k_\tau}\|
 + L_\lambda L_c\alpha_{T,k_\tau}^2\|g_{T,k_\tau}\|^2.
\end{aligned}
\eeqn
Taking \req{switch} into account
gives $\|c_{k_\tau}\|\le \beta \alpha_{T,k_\tau}\|g_{T,k_\tau}\|$, and
we obtain 
\[
\Delta_\lambda
\le \left(\beta L_\lambda  + L_\lambda L_c\right)\alpha_{T,k_\tau}^2\|g_{T,k_\tau}\|^2.
\]
Thus, summing $\Delta_x$ and $\Delta_\lambda$,
\[
\begin{aligned}
\psi(x_{k_\tau+1}) - \psi(x_{k_\tau})
\leq & -\alpha_{T,k_\tau}\|g_{T,k_\tau}\|^2
+\frac{1}{2}(L_L+\rho L_c)\alpha_{T,k_\tau}^2\|g_{T,k_\tau}\|^2\\
&+\left(\beta L_\lambda  + L_\lambda L_c\right)\alpha_{T,k_\tau}^2\|g_{T,k_\tau}\|^2
\end{aligned}
\]
and \req{effectT} follows.
} \epr

\noindent
Bound \req{effectT} quantifies the effect of  tangential steps
on the Lyapunov function, and its right-hand side involves a
first-order (descent) term and a second-order perturbation term.
We now derive crucial bounds on these terms.

\llem{lem:adagrad}{
If we denote
\[
\Gamma_{k_{\tau+1}} =\Gamma_{k_{\tau}}+\|g_{T,k_{\tau}}\|^2
\tim{ and }
\alpha_{T,k_\tau} =\frac{\eta}{\sqrt{\Gamma_{k_{\tau+1}}+\varsigma}},
\]
then, for all $0 \le \tau_0\le \tau_1$,
\beqn{eq:AG1}
\sum_{\tau=\tau_0}^{\tau_1}\alpha_{T,k_\tau}\,\|g_{T,k_\tau}\|^2
\ge \eta\,\left(\sqrt{\Gamma_{k_{\tau_1+1}}+\varsigma}-\sqrt{\Gamma_{k_{\tau_0}}+\varsigma}\right),
\eeqn
\beqn{eq:AG2}
\sum_{\tau=\tau_0}^{\tau_1}\alpha_{T,k_\tau}^2\,\|g_{T,k_\tau}\|^2
\le
\eta^2\,\log\left(\frac{\Gamma_{k_{\tau_1+1}}+\varsigma}{\Gamma_{k_{\tau_0}}+\varsigma}\right).
\eeqn
}

\proof{
Let $w_{k_\tau+1} =\sqrt{\Gamma_{k_{\tau+1}}+\varsigma}$.
The definition \req{alphaT-def} implies
\[
\begin{aligned}
\sum_{\tau=\tau_0}^{\tau_1} \alpha_{T,k_\tau} \,\|g_{T,k_\tau}\|^2
& =    \eta \sum_{\tau=\tau_0}^{\tau_1}
\frac{\|g_{T,k_\tau}\|^2}{\sqrt{\Gamma_{k_{\tau+1}}+\varsigma}}\\
&> \eta \sum_{\tau=\tau_0}^{\tau_1} \frac{\|g_{T,k_\tau}\|^2}{w_{k_{\tau+1}}+w_{k_\tau}}\\
&=    \eta \sum_{\tau=\tau_0}^{\tau_1} \frac{w_{k_{\tau+1}}^2-w_{k_{\tau}}^2}{w_{k_{\tau+1}}+w_{k_{\tau}}}\\
&=    \eta \sum_{\tau=\tau_0}^{\tau_1} (w_{k_{\tau+1}}-w_{k_{\tau}})\\
&=    \eta  (w_{k_{\tau_1+1}}-w_{k_{\tau_0}}),
\end{aligned}
\]
which gives \req{eq:AG1}.  Using the concavity and the increasing
nature of the logarithm, we also
deduce from \req{alphaT-def} that
\[
\sum_{\tau=\tau_0}^{\tau_1}\alpha_{T,k_\tau}^2\|g_{T,k_\tau}\|^2
= \eta^2 \!\sum_{\tau=\tau_0}^{\tau_1}\! \frac{\|g_{T,k_\tau}\|^2}{\Gamma_{k_{\tau+1}}+\varsigma}
= \eta^2 \!\sum_{\tau=\tau_0}^{\tau_1}\! \frac{\Gamma_{k_{\tau+1}}-\Gamma_{k_{\tau}}}{\Gamma_{k_{\tau+1}}+\varsigma}
\leq \eta^2 \!\sum_{\tau=\tau_0}^{\tau_1}\!
\log(\Gamma_{k_{\tau+1}}+\varsigma)-\log(\Gamma_{k_{\tau}}+\varsigma),
\]
giving \req{eq:AG2}.
}  \epr

\noindent
Note that the inequalities \req{eq:AG1} and \req{eq:AG2} remain valid if there are
no tangential iterations, since $\Gamma_{k_{\tau_1+1}} = 0$ in this
case.

\subsection{Telescoping sum}

We now combine the effects of tangential and normal steps to
accumulate in a telescoping sum the relevant decreases on $\psi(x)$ across
all iterations. 

\llem{telescoping}{
Suppose A0--A6 hold. Then
\beqn{full-telescopic}
\sqrt{\Gamma_{k_{\tau_1+1}}+\varsigma} + \sum_{\nu=\nu_0}^{\nu_1}\|c_{k_\nu}\|
\leq \kap{gap} + \kap{tan}\log\left(1+\frac{\Gamma_{k_{\tau_1+1}}}{\varsigma}\right),
\eeqn
where
\[
\kap{gap} = \frac{1}{\eta}\left(\psi(x_0) +
\frac{\kappa_J\kappa_g\kappa_c}{\sigma_0^2} +\rho\kappa_c - \flow\right)+\sqrt{\varsigma}.
\]
}

\proof{
Combining \req{effectT}, \req{eq:AG1} and \req{eq:AG2}, we obtain
\beqn{eq:sumTtg}
\begin{aligned}
\sum_{\tau=\tau_0}^{\tau_1}\left(\psi(x_{k_\tau+1})-\psi(x_{k_\tau})\right)
&\le - \eta\,\sqrt{\Gamma_{k_{\tau_1+1}}+\varsigma}+ \eta\,\sqrt{\Gamma_{k_{\tau_0}}+\varsigma}
 + \eta^2\kap{tan}\log\left(\frac{\Gamma_{k_{\tau_1+1}}+\varsigma}{\Gamma_{k_{\tau_0}}+\varsigma}\right).
\end{aligned}
\eeqn
Then, if
$\min[k_{\nu_0},k_{\tau_0}] = 0$ and $\max[k_{\nu_1},k_{\tau_1}]= k$,
we have $w_{k_{\tau_0}} = \varsigma$ and $\Gamma_{k_{\tau_0}}=0$.
The bounds  \req{effectN} and \req{eq:sumTtg} therefore give
\beqn{eq:sumTtg2}
\begin{aligned}
\psi(x_{k+1}) &- \psi(x_0)\\
&=\sum_{\tau=\tau_0}^{\tau_1}\big(\psi(x_{k_\tau+1})-\psi(x_{k_\tau})\big)
+ \sum_{\nu=\nu_0}^{\nu_1}\big(\psi(x_{k_\nu+1})-\psi(x_{k_\nu})\big) \\
& \le - \eta\,\sqrt{\Gamma_{k_{\tau_1+1}}+\varsigma}+\eta\sqrt{\varsigma}
      - \eta \sum_{\nu=\nu_0}^{\nu_1}\|c_{k_\nu}\|
      + \eta^2\kap{tan}\log\left(\frac{\Gamma_{k_{\tau_1+1}}+\varsigma}{\Gamma_{k_{\tau_0}}+\varsigma}\right)\\
& \le - \eta\,\sqrt{\Gamma_{k_{\tau_1+1}}+\varsigma}+\eta\sqrt{\varsigma}
      - \eta \sum_{\nu=\nu_0}^{\nu_1}\|c_{k_\nu}\|
      + \eta\kap{tan}\log\left(1+\frac{\Gamma_{k_{\tau_1+1}}}{\varsigma}\right),
\end{aligned}
\eeqn
where we used the fact that $\eta \le 1$.
But \req{eq:ls-mult}, A2, A3, and A4 ensure that
\[
\|\widehat\lambda(x)\|\leq \frac{\kappa_J}{\sigma_0^2}\|g(x)\|
\leq \frac{\kappa_J\kappa_g}{\sigma_0^2}.
\]
Hence, using \req{lyap-def}, A1, A2 and Lemma~\ref{c-bounded}, we have
\[
\psi(x_{k+1}) -\psi(x_0) \geq \flow-
\frac{\kappa_J\kappa_g\kappa_c}{\sigma_0^2}-\rho\kappa_c-\psi(x_0)
\eqdef -\eta(\kap{gap}-\sqrt{\varsigma}),
\]
so that \req{eq:sumTtg2} imply \req{full-telescopic}.
} \epr

\noindent
This result is central in our analysis as it provides consistent upper
bounds on a quantity related to  the optimality
measure $\|g_{T,k}\|+\|c_k\|$, which we now exploit separately.

\subsection{Tangential complexity}

Our next step is to use \req{full-telescopic} to derive a global rate
of convergence over tangential steps.  We start by proving a useful
technical result.

\llem{tech-tan}{
  Suppose that $at \le b + c\,\log(t)$  for $t\geq 1$.  Then
  \[
  t \le \max\left[e^{b/c}, \frac{4c^2}{a^2}\right].
  \]
}

\proof{ Suppose that $t> e^{b/c}$.  Then $b\le c\,\log(t)$ and thus
  (see \cite[eq. (14)]{Tops07})
\[
at \le 2c\,\log(t) \le 2c\,\log(1+t) \le \frac{2ct}{\sqrt{1+t}}.
\]
Hence $a\sqrt{1+t} \le 2c$, 
which is to say $t \le \left(2c/a\right)^2-1$, yielding
the desired bound.
}

\noindent
We now consider the rate of convergence for tangential steps proper.

\llem{tangent-complexity}{
Suppose A0--A6 hold.  Then
\beqn{Gamma-bound}
\sqrt{\Gamma_{k_{\tau_1+1}}+\varsigma} \le \kappa_T
\eqdef
\sqrt{\varsigma}\max\left[e^{\frac{\kap{gap}}{2\kap{tan}}},
  \frac{16\kap{tan}^2}{\varsigma}\right]
\eeqn
and
\beqn{Tcomp}
\sum_{\tau=\tau_0}^{\tau_1}\|g_{T,k_\tau}\|+\|c_{k_\tau}\|
\leq \kappa_T\sqrt{\tau_1+1}\left(1+\frac{\beta\eta}{\sqrt{\varsigma}}\right).
\eeqn
}

\proof{
The bound \req{full-telescopic} implies
\[
\sqrt{\Gamma_{k_{\tau_1+1}}+\varsigma}
\leq \kap{gap} +\kap{tan}\log\left(1+\frac{\Gamma_{k_{\tau_1+1}}}{\varsigma}\right).
\]
and thus
\[
\sqrt{\varsigma}\sqrt{1+\frac{\Gamma_{k_{\tau_1+1}}}{\varsigma}}
\leq \kap{gap} +\kap{tan}\log\left(1+\frac{\Gamma_{k_{\tau_1+1}}}{\varsigma}\right)
= \kap{gap} +2\kap{tan}\log\left(\sqrt{1+\frac{\Gamma_{k_{\tau_1+1}}}{\varsigma}}\right).
\]
Using Lemma~\ref{tech-tan}, we then obtain
\[
\sqrt{\Gamma_{k_{\tau_1+1}}+\varsigma}
= \sqrt{\varsigma}\sqrt{1+\frac{\Gamma_{k_{\tau_1+1}}}{\varsigma}}
\le
\sqrt{\varsigma}\max\left[e^{\frac{\kap{gap}}{2\kap{tan}}},
  \frac{16\kap{tan}^2}{\varsigma}\right]
=\kappa_T.
\]
This is \req{Gamma-bound}.
We may now invoke the inequality
\[
\sum_{j=0}^k a_j \le \sqrt{k+1}\sqrt{\sum_{j=0}^k a_j^2}
\]
for nonnegative $\{a_j\}_{j=0}^k$ to deduce from \req{Gamma-def} and \req{Gamma-bound} that
\beqn{tantan}
\sum_{\tau=\tau_0}^{\tau_1} \|g_{T,k_\tau}\|
\le \sqrt{\tau_1+1}\,\sqrt{\sum_{\tau=\tau_0}^{\tau_1}\|g_{T,k_\tau}\|^2}
= \sqrt{\tau_1+1}\,\sqrt{\Gamma_{k_{\tau_1+1}}}
\le \sqrt{\tau_1+1}\,\kappa_T.
\eeqn
Using the switching condition \req{switch} and the fact that $
\alpha_{T,k_\tau} \leq \eta/\sqrt{\varsigma}$, we also deduce that
\[
\sum_{\tau=\tau_0}^{\tau_1}\|c_{k_\tau}\|
\le \sum_{\tau=\tau_0}^{\tau_1} \beta \alpha_{T,k_\tau}\|g_{T,k_\tau}\|
\le \frac{\beta\eta}{\sqrt{\varsigma}} \sum_{\tau=\tau_0}^{\tau_1} \|g_{T,k_\tau}\|
\le \frac{\beta\eta\kappa_T\sqrt{\tau_1+1}}{\sqrt{\varsigma}}.
\]
Summing this bound with \req{tantan} then gives \req{Tcomp}.
} \epr

\noindent
A marginally tighter bound on $\Gamma_{k_{\tau_1}}$ can
obtained, as in \cite{BellGratMoriToin25}, by using the Lambert $W_{-1}$
function (see \cite{Corletal96}) instead of applying
Lemma~\ref{tech-tan}, but we have chosen the latter for
simplicity of exposition.

\subsection{Normal complexity}

The analysis of the complexity of normal steps also uses the switching
condition, but in the other direction.

\llem{normal-complexity}{
Suppose A0--A6 hold. Then
\beqn{Ncomp}
\sum_{\nu=\nu_0}^{\nu_1}\|g_{T,k_\nu}\|+\|c_{k_\nu}\|
< \kappa_N\,\left(1+\frac{\kappa_T+\kappa_g}{\beta\eta}\right),
\eeqn
where
\beqn{kappaN-def}
\kappa_N = \kap{gap}+2\kap{tan}\log\left(\frac{\kappa_T}{\sqrt{\varsigma}}\right).
\eeqn
}

\proof{
The bound \req{full-telescopic}  ensures that
\beqn{nc1-sto}
\sum_{\nu=\nu_0}^{\nu_1}\|c_{k_\nu}\|
\leq \kap{gap}+\kap{tan}\log\left(1+\frac{\Gamma_{k_{\tau_1+1}}}{\varsigma}\right),
\eeqn
where $k_{\tau_1}$ is the index of the last tangential iteration before
$k_{\nu_1}$. As in Lemma~\ref{tangent-complexity}, we may
now use the
 monotonicity of the logarithm to obtain
\beqn{kappaN-bound}
\sum_{\nu=\nu_0}^{\nu_1}\|c_{k_\nu}\|
\le \kap{gap}+2\kap{tan}\log\left(\sqrt{1+\frac{\Gamma_{k_{\tau_1}}}{\varsigma}}\right)
\le\kap{gap}+2\kap{tan}\log\left(\frac{\kappa_T}{\sqrt{\varsigma}}\right)
= \kappa_N.
\eeqn
Using the switching condition
\req{switch}, we obtain, for $\nu \in \iibe{\nu_0}{\nu_1}$,
\beqn{nc01}
\|c_{k_\nu}\| \ge \beta\,\alpha_{T,k_\nu}\|g_{T,k_\nu}\|.
\eeqn
Observe now that, using \req{Gamma-bound}, A2 and $\sqrt{a+b}\le\sqrt{a}+\sqrt{b}$, 
\[
\alpha_{T,k_\nu}
= \frac{\eta}{\sqrt{\Gamma_{k_\nu}+\|g_{T,k}\|^2+\varsigma}}
\ge \frac{\eta}{\sqrt{\Gamma_{k_\nu}+\varsigma}+\kappa_g}
\ge \frac{\eta}{\kappa_T+\kappa_g}.
\]
Thus, \req{nc01} yields
\beqn{nc1b}
\|c_{k_\nu}\|\
\ge \frac{\beta\eta}{\kappa_T+\kappa_g}\|g_{T,k_\nu}\|
\eeqn
for $\nu \in \iibe{\nu_0}{\nu_1}$.
With \req{kappaN-bound}, this implies
\[
\sum_{\nu=\nu_0}^{\nu_1}\|g_{T,k_\nu}\|
\le \frac{\kappa_T+\kappa_g}{\beta\eta}\sum_{\nu=\nu_0}^{\nu_1}\|c_{k_\nu}\|
\le \frac{\kappa_N\,(\kappa_T+\kappa_g)}{\beta\eta}.
\]
Summing this bound with \req{kappaN-bound} gives \req{Ncomp}. 
} \epr

\subsection{Combined complexity}

The combined complexity may now be derived by assembling the above results.

\lthm{final-complexity}{
  Suppose A0-A6 hold. Then
\beqn{the-complexity}
\frac{1}{k+1}\sum_{j=0}^k\Big(\|g_{T,j}\|+\|c_j\|\Big)
\leq \frac{\kap{ADSW,1}}{\sqrt{k+1}} + \frac{\kap{ADSW,2}}{k+1}
= \calO\left(\frac{1}{\sqrt{k+1}}\right),
\eeqn
where
\[
\kap{ADSW,1}= \kappa_T\left(1+\frac{\beta\eta}{\sqrt{\varsigma}}\right)
\tim{ and }
\kap{ADSW,2}= \kappa_N\,\left(1+\frac{\kappa_T+\kappa_g}{\beta\eta}\right).
\]
}

\proof{
Consider iterations from $0$ to $k$ of both types (tangential
and normal) by setting $\min[k_{\nu_0},k_{\tau_0}] = 0$ and
$\max[k_{\nu_1},k_{\tau_1}]= k$ (as in Lemma~\ref{telescoping}). We then
obtain, by combining \req{Tcomp} and \req{Ncomp},
\[
\begin{aligned}
\sum_{j=0}^k\big(\|g_{T,j}\|+\|c_j\|\Big)&
= \sum_{\tau=\tau_0}^{\tau_1}\Big(\|g_{T,k_\tau}\|+\|c_{k_\tau}\|\Big)
+\sum_{\nu=\nu_0}^{\nu_1}\Big(\|g_{T,k_\nu}\|+\|c_{k_\nu}\|\Big)\\
&\le \kappa_T\sqrt{k+1}\left(1+\frac{\beta\eta}{\sqrt{\varsigma}}\right)
+\kappa_N\,\left(1+\frac{\kappa_T+\kappa_g}{\beta\eta}\right),
\end{aligned}
\]
where we used the inequalities  $\tau_1 \le k_{\tau_1}\leq k$ and
$k_{\nu_1}\leq k$. Bound \req{the-complexity} is finally obtained
by dividing both sides by $k+1$.
}\epr

\numsection{Numerical illustration}\label{sec:numerics}

To illustrate the behaviour of \al{ADSWITCH}, we coded the algorithm in
Matlab and applied it on a set of small problems from the
{\sf CUTEst} collection \cite{GoulOrbaToin15b} as supplied in Matlab
by S2MPJ \cite{GratToin24}.  Our implementation uses the Newton normal
step \req{facto-Newton-normal} and the constants
\[
\beta = 0.01,
\ms
\eta = 1,
\ms
\theta = 1000,
\tim{and}
\delta = \varsigma= 10^{-5}.
\]
For a given $\epsilon\in (0,1)$, a run on a given test problem is deemed successful at iteration $k$ 
if
\beqn{convg}
\max[\|g_{T,k}\|,\|c_k\|] \le \epsilon  \tim{(convergence)}
\eeqn
or
\beqn{infeas}
\|J_k^Tc_k\|\le \epsilon \tim{but} \|c_k\|>\epsilon \tim{(infeasible critical point)}
\eeqn
or if a sufficiently optimal function value was found in the sense that
\beqn{accept}
\|c_k\|\le \epsilon
\tim{and}
\left\{\begin{array}{ll}
|f(x_k)| \le |f_*| + 10^{-7}   &\tim{if} |f_*| < 10^{-7} \\
|f(x_k)-f_*| \le 10^{-7} |f_*| &\tim{if} |f_*| \ge 10^{-7}, \\
\end{array}\right.
\eeqn
where $f_*$ is the best known feasible function value for the problem.
Optimization was also stopped after a maximum of 100000
iterations. 

We first ran all test problems with accurate gradient values
and $\epsilon = 10^{-5}$. The complete results are reported in
the appendix, and may be broadly summarized as follows.
\begin{enumerate}
\item The algorithm appears to converge as predicted.  Despite the
  fact that our current analysis does not cover the case where the
  Jacobian may lose rank, the algorithm did solve a number of 
  instances where this occurs, sometimes finding an infeasible
  critical point of the constraints violation\footnote{This occurs,
  for instance, at the starting point of HS61.}. 
\item Its performance and reliability is (as could be anticipated)
  dominated by that of AdaGrad for the tangential step.  Since this
  method is purely first-order, it often performs well, but may fail
  to solve ill-conditioned problems in a reasonable number of
  iterations. This is also the case for \al{ADSWITCH}. Over the 71
  problems in our test set, it solves 44 (62\%) of them within 750
  iterations, 58 (81\%) of them within 100000 iterations, and fails on
  13. While reliability is not a strong point of AdaGrad and
  \al{ADSWITCH} in the deterministic case, the situation is very
  different when noise is present in the objective function's gradient, as we show below.
\end{enumerate}

\begin{figure}[htbp]
\includegraphics[width=0.48\textwidth]{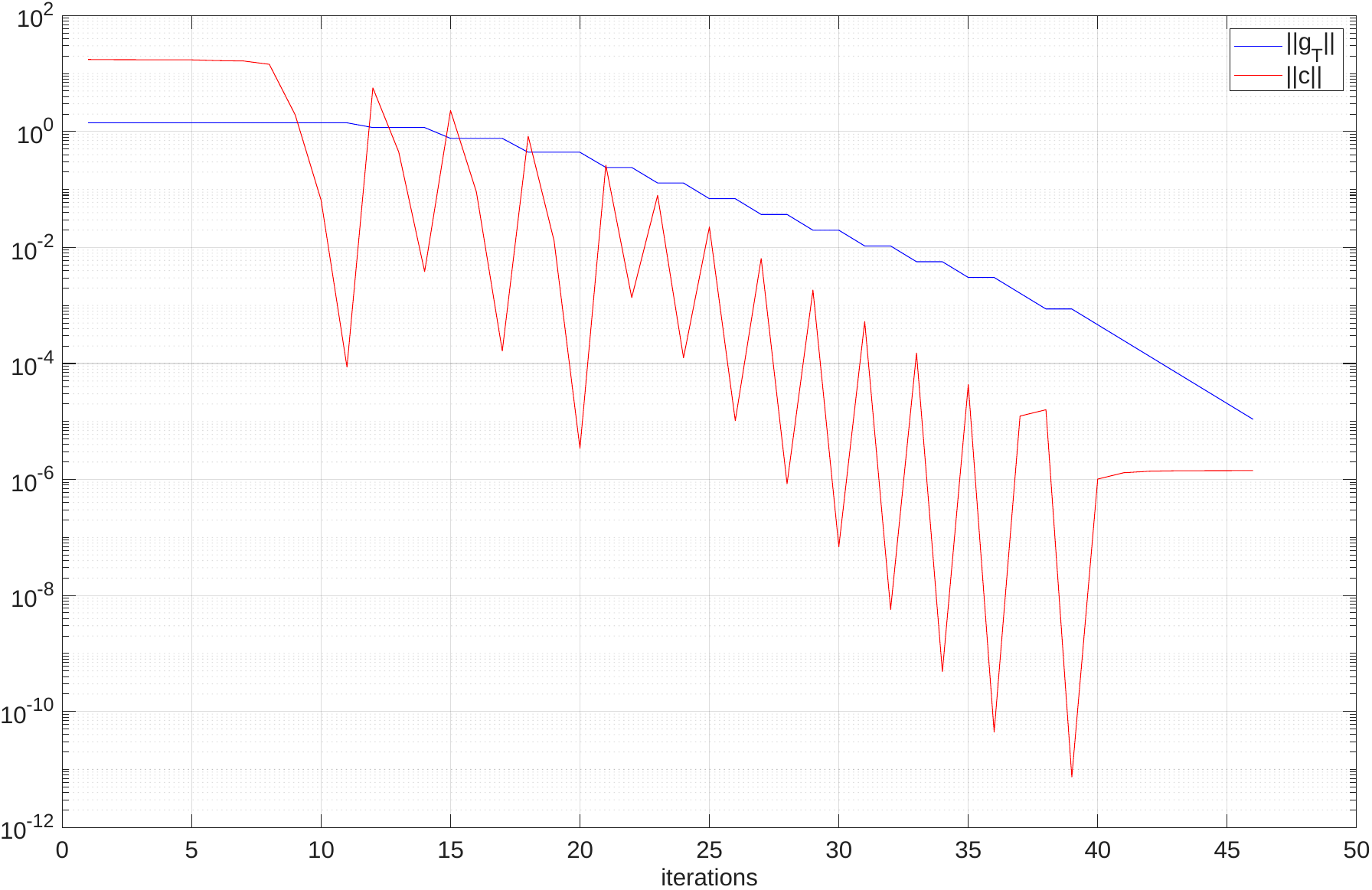}
\includegraphics[width=0.48\textwidth]{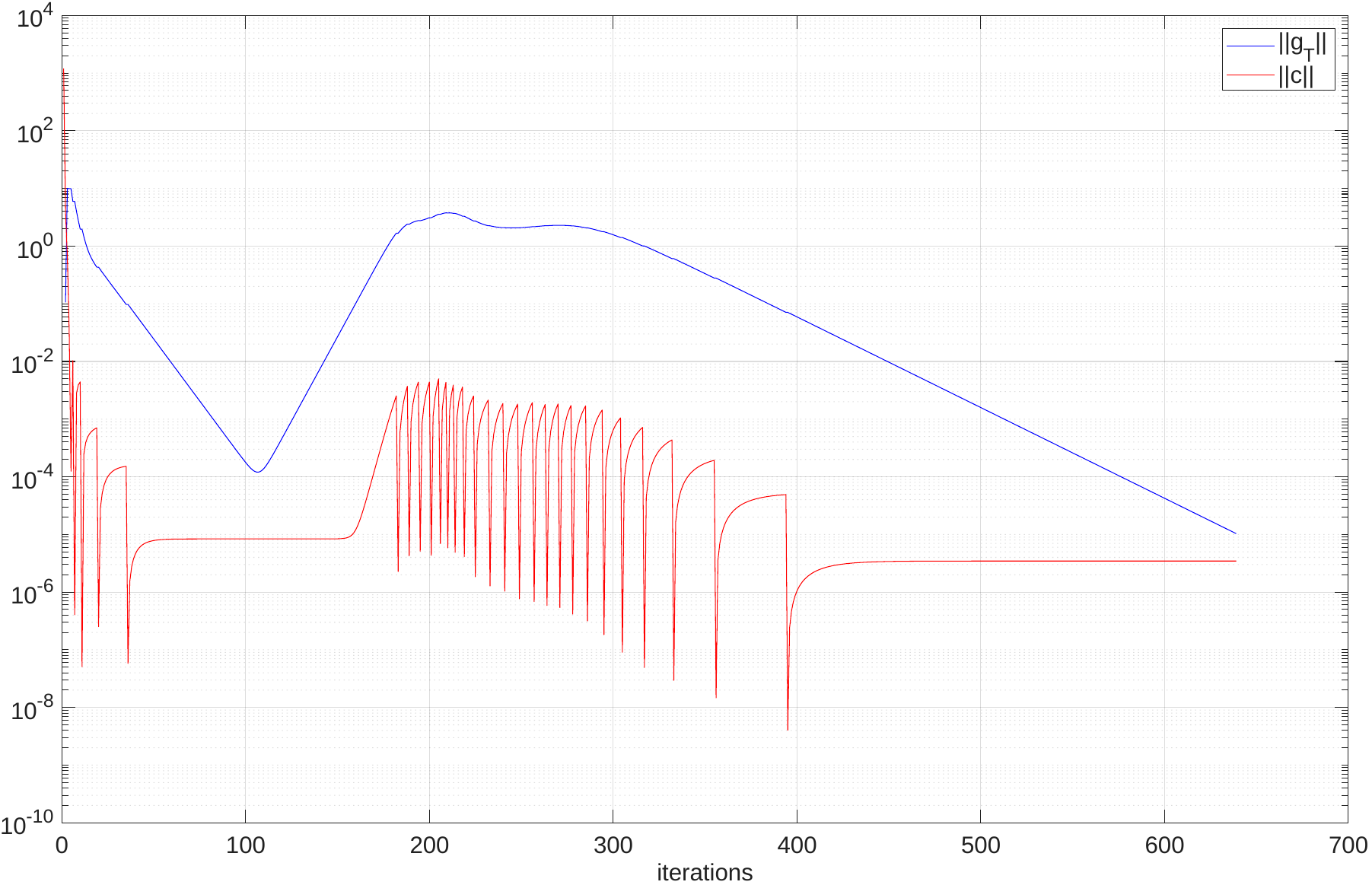}
\caption{\label{oks-1-2}Evolution of the
  optimality measure (left: {\tt BYRDSPHR}, right: {\tt ORTHREGA})}
\end{figure}

\begin{figure}[htbp]
\includegraphics[width=0.48\textwidth]{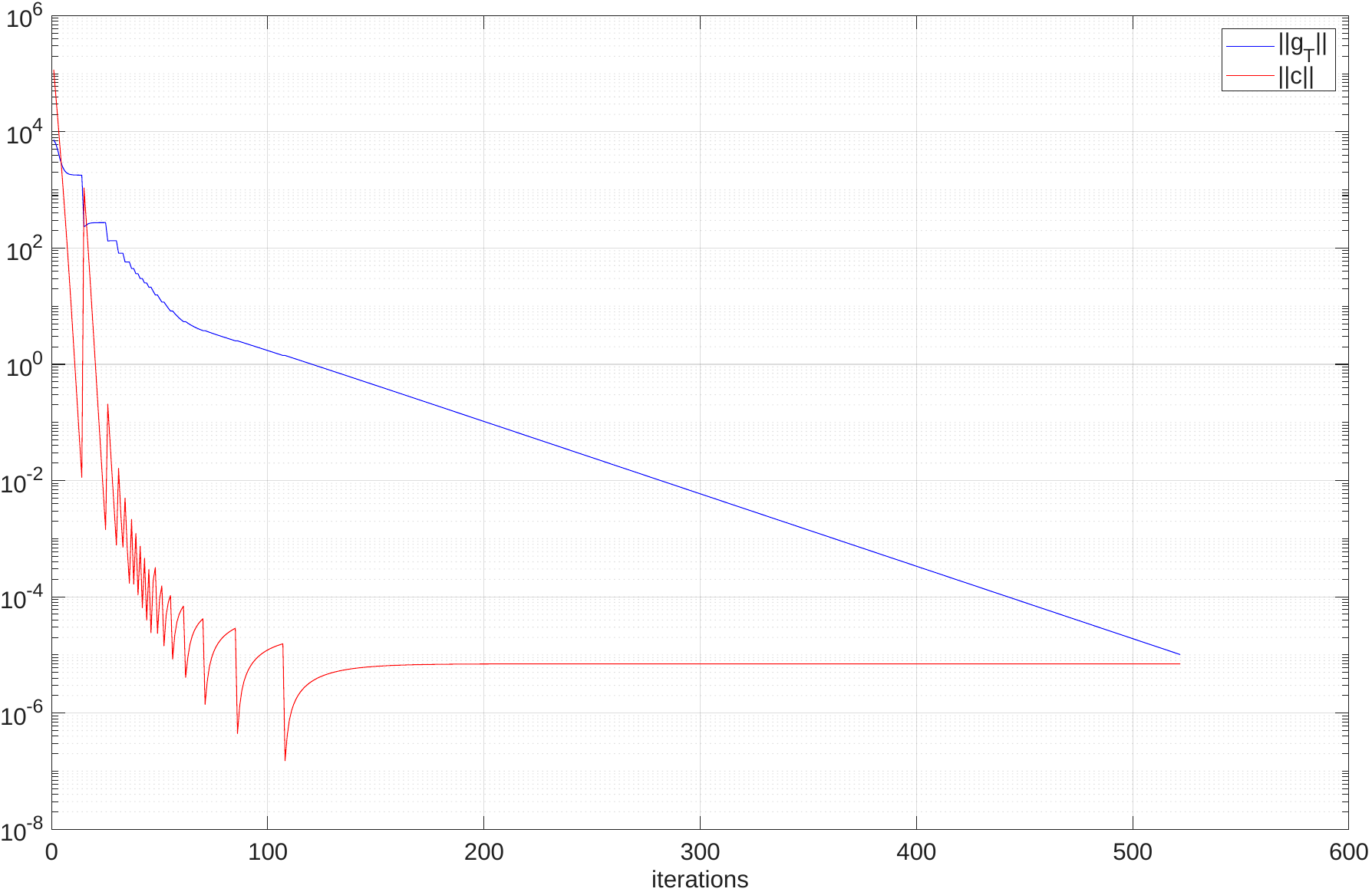}
\includegraphics[width=0.48\textwidth]{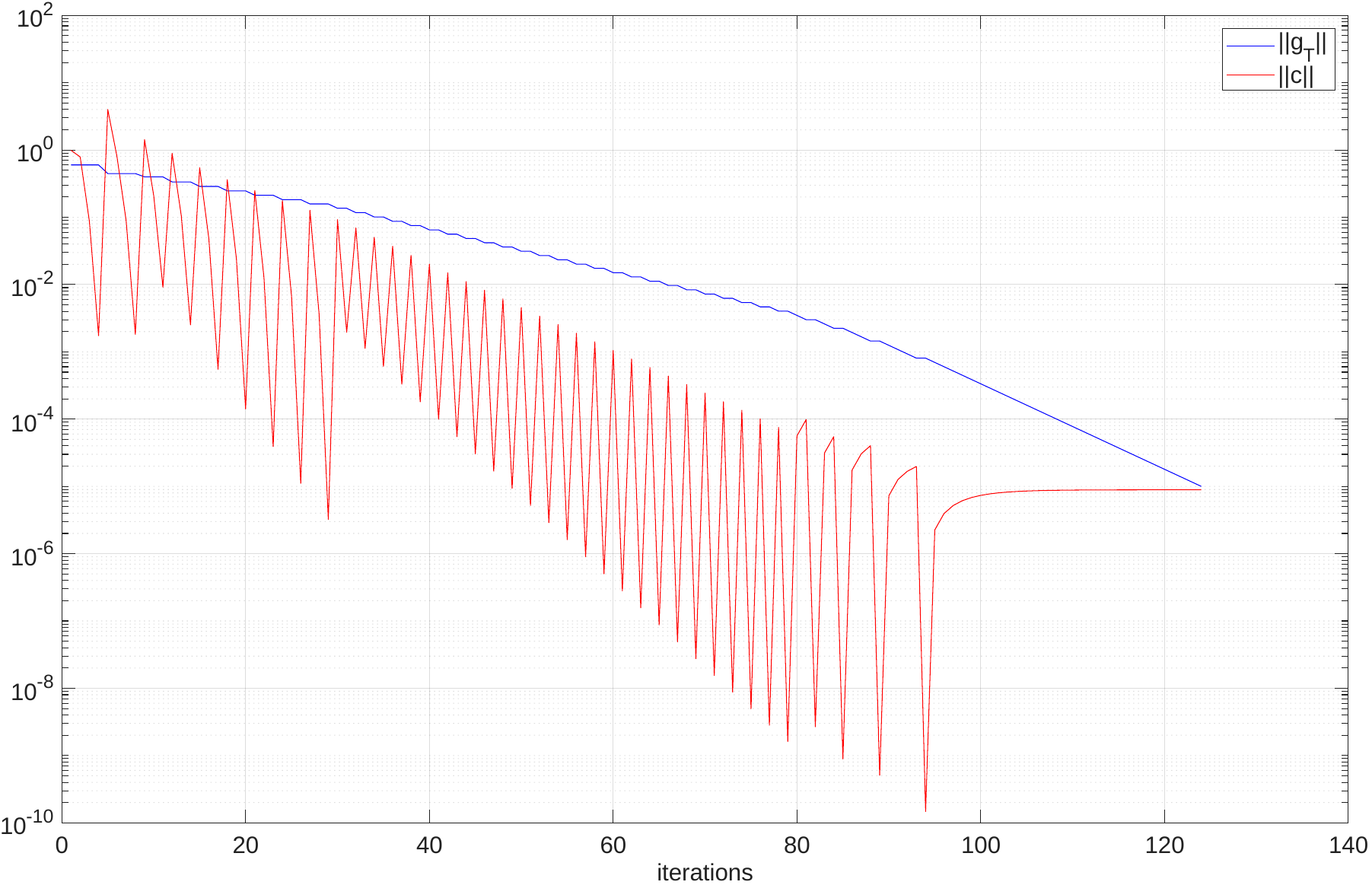}
\caption{\label{oks-3-4}Evolution of the
  optimality measure (left: {\tt LUKVLE2}, right: {\tt BT1})}
\end{figure}

\begin{figure}[htbp]
\begin{minipage}{0.45\textwidth}
   \centerline{\includegraphics[width=0.9\textwidth]{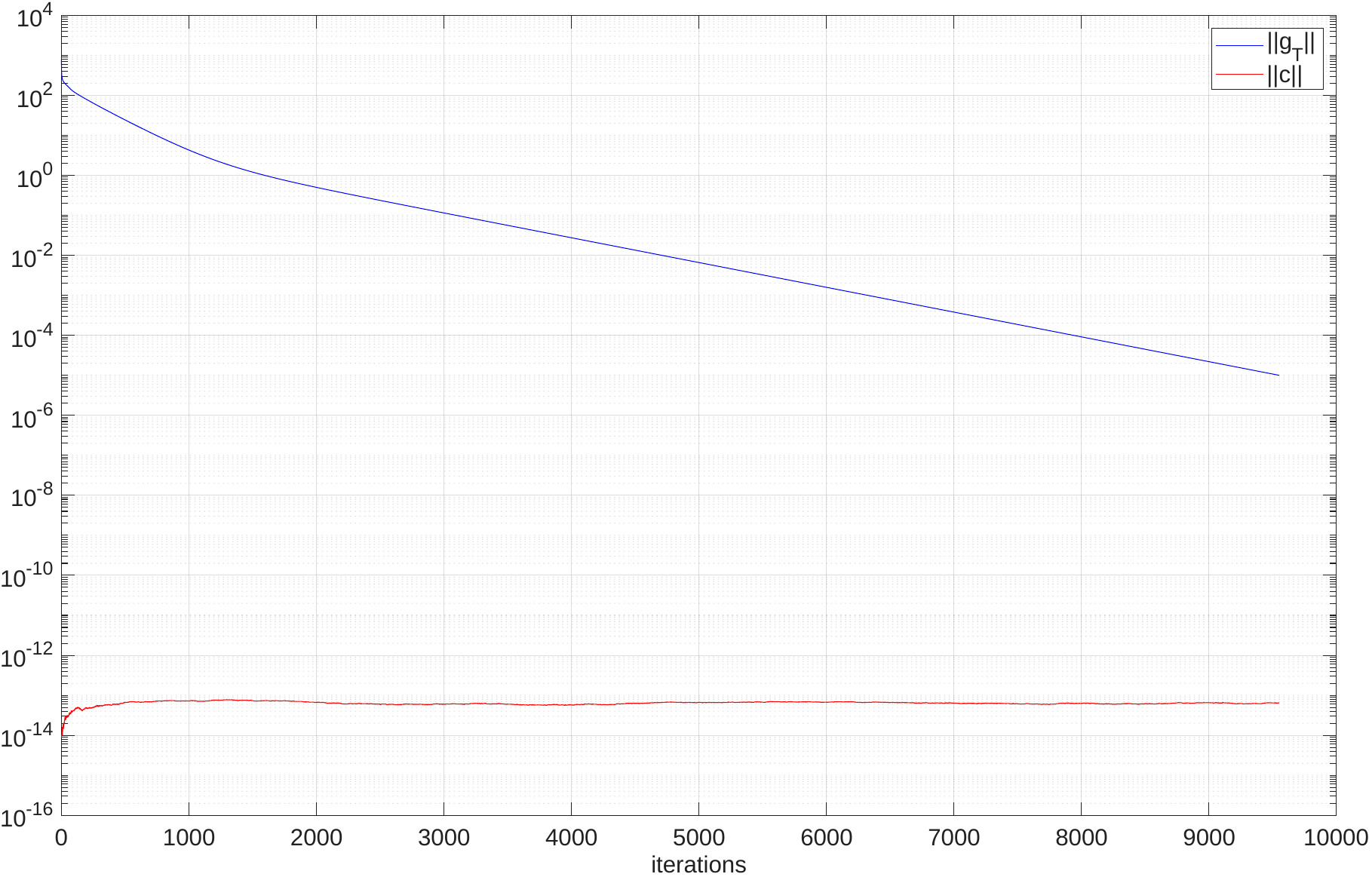}}
   \caption{\label{bad}Evolution  of the
            optimality measure ({\tt HS50})}
\end{minipage}
\hspace*{0.08\textwidth}  
\begin{minipage}{0.45\textwidth}
    \vspace*{1cm}
    \begin{tabular}{|c|c|c|}
      \hline
      {\footnotesize relative}     &  {\footnotesize number of} &
      {\footnotesize number of} \\
       {\footnotesize noise level} &  {\footnotesize total failures}
       &  {\footnotesize total successes} \\
      \hline
      5\%  &  6 & 61 \\
      15\% &  7 & 53 \\
      25\% &  7 & 53 \\
      50\% & 12 & 51 \\
      \hline
    \end{tabular}
    \vspace*{1cm}
    \caption{\label{tab:reliab}Reliability of \al{ADSWITCH} 
    in the presence of relative Gaussian noise in the objective
    functions's gradient.}
\end{minipage}
\end{figure}

We illustrate a few cases of satisfactory convergence in
Figures~\ref{oks-1-2} and \ref{oks-3-4}.  One can
clearly see the difference in speed of convergence between tangential
steps (AdaGrad-like) and normal steps (Newton), and thus that the
overall performance is dominated by that of the first-order method
defining the tangential step. This is even clearer when convergence is
too slow, as shown in Figure~\ref{bad} where constraint violation
remains very small\footnote{{\tt HS50} has linear
equality constraints only.}, while very slow
convergence of the projected gradient results in a large number of iterations.
The right panel of Figure~\ref{oks-1-2} ({\tt ORTHREGA} with $n=133$ and
$m=64$) is also interesting because it shows the
switching condition \req{switch} in action: the algorithm first
approaches the neighbourhood of an infeasible critical point of the constraints violation
until the switching rule allows the tangential step to take over,
causing the iterates to escape this neighbourhood.


\begin{figure}[htbp]
\includegraphics[width=0.48\textwidth]{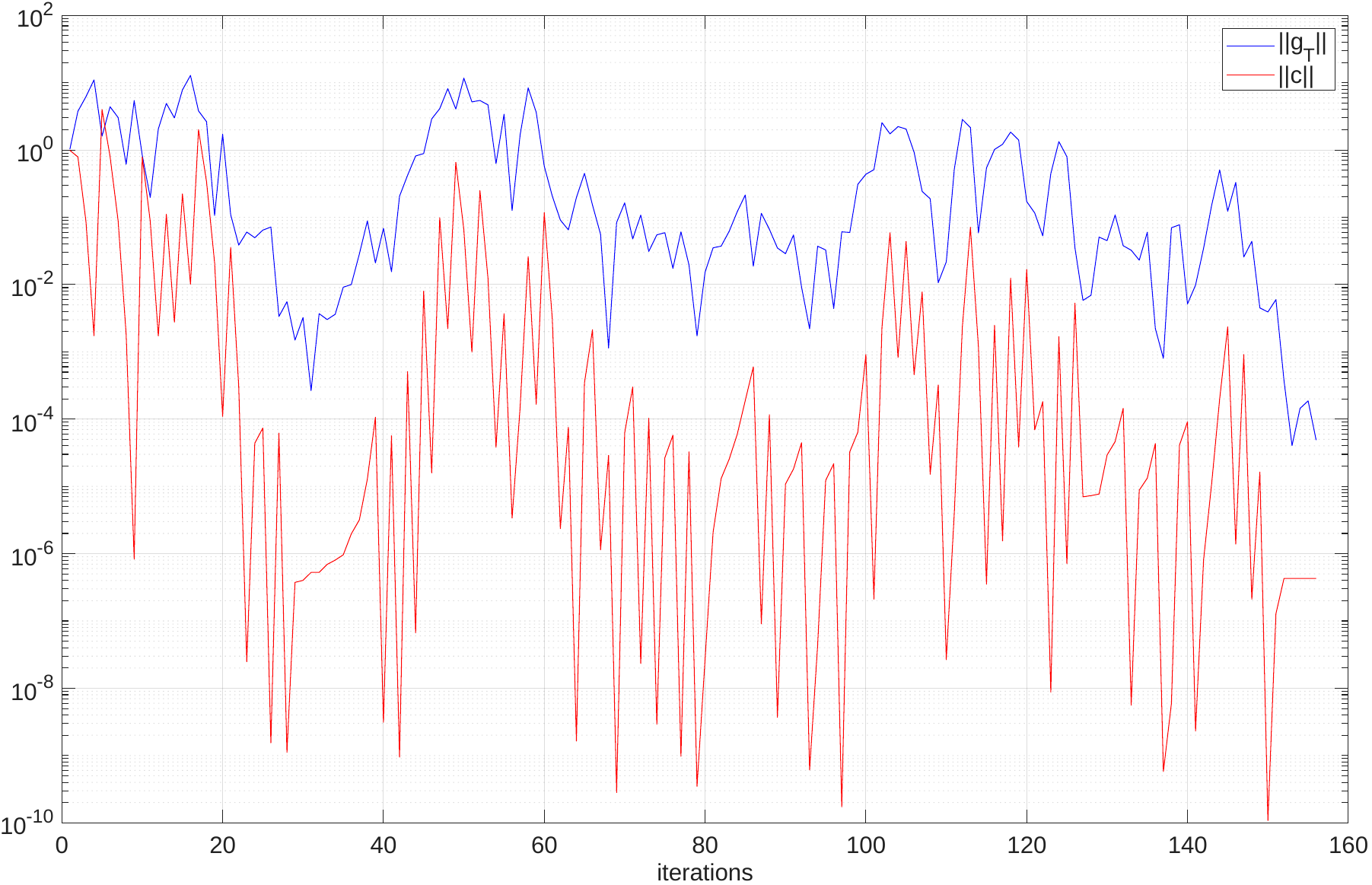}
\includegraphics[width=0.48\textwidth]{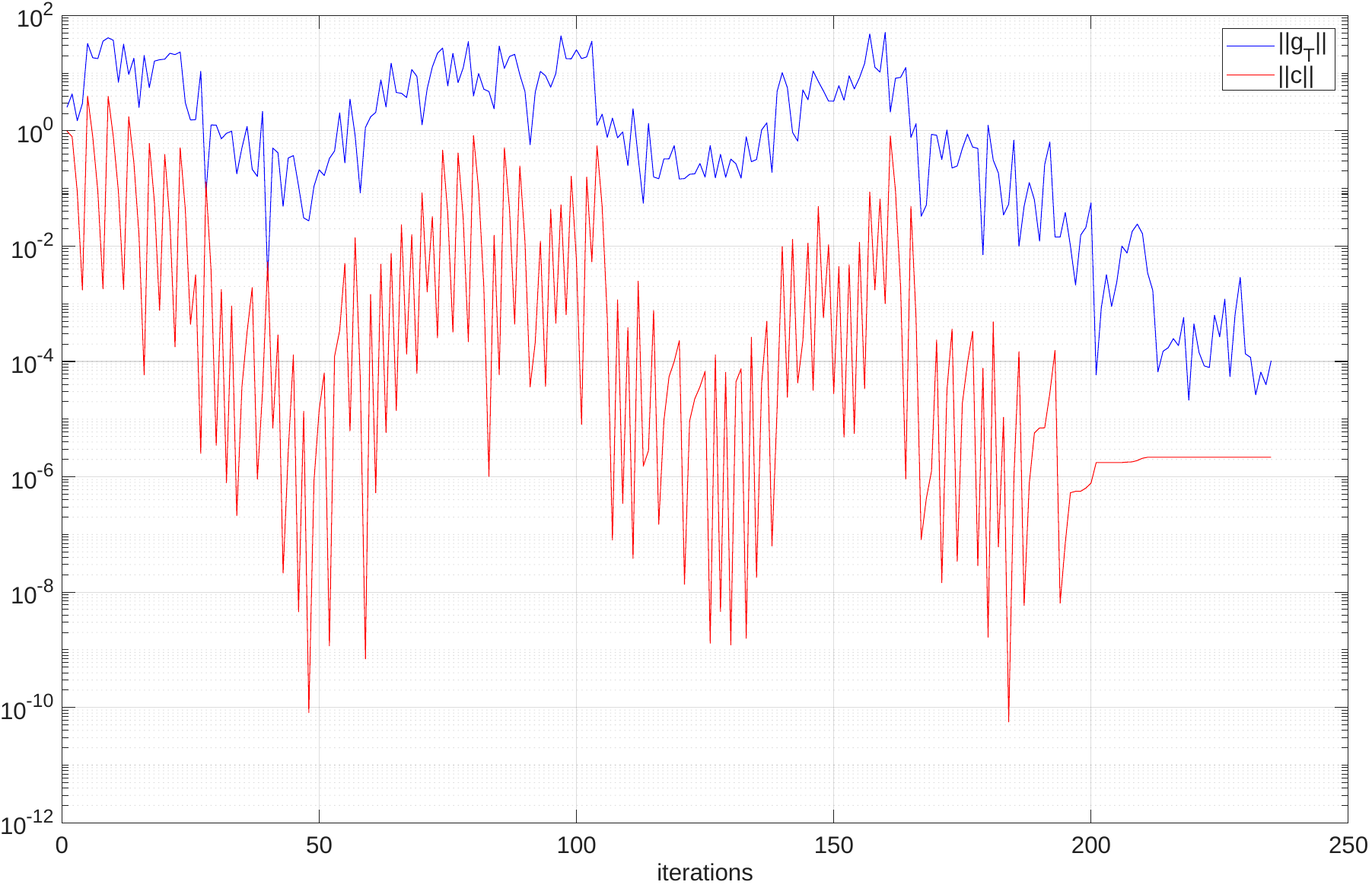}\\
\includegraphics[width=0.48\textwidth]{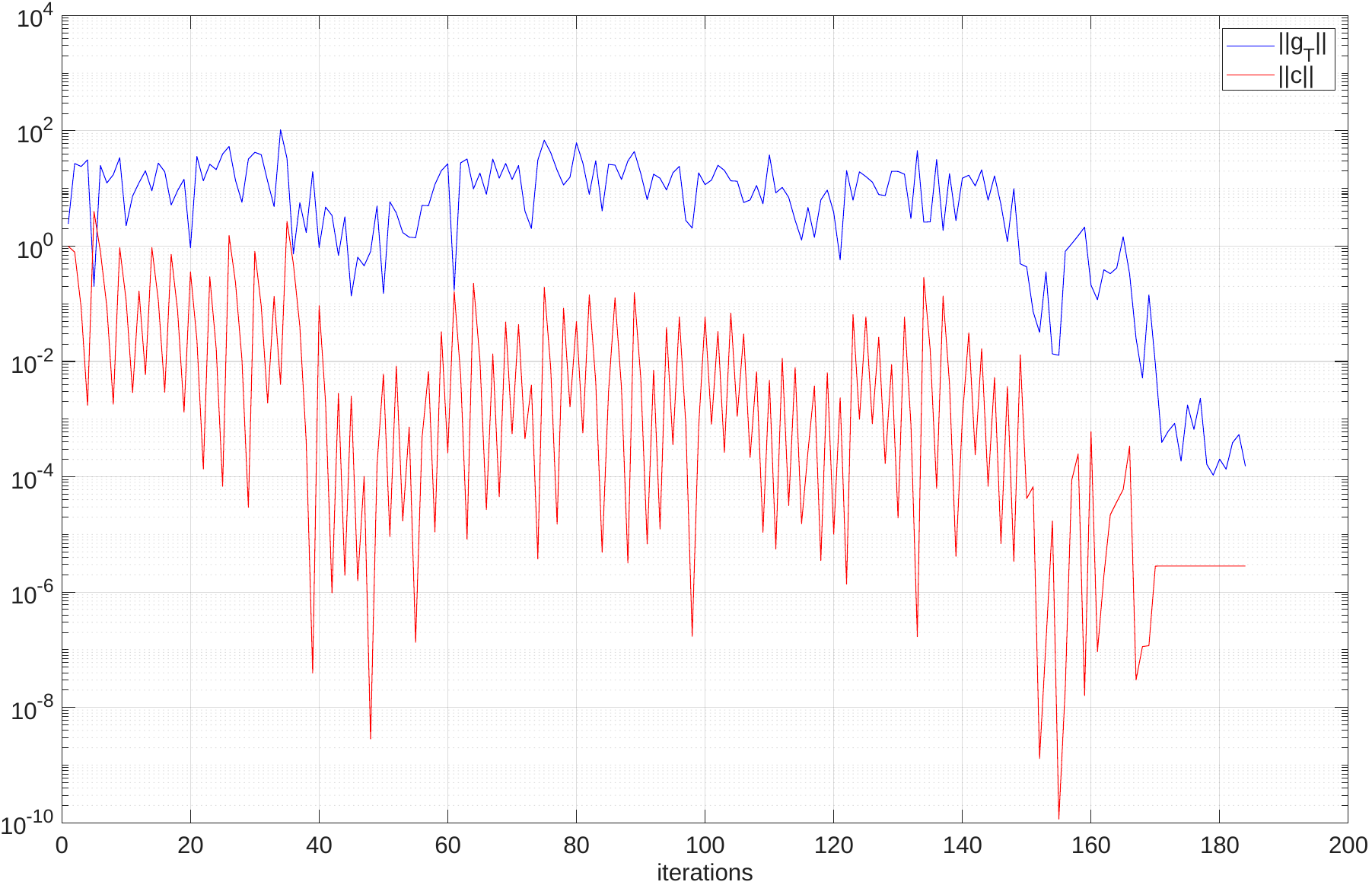}
\includegraphics[width=0.48\textwidth]{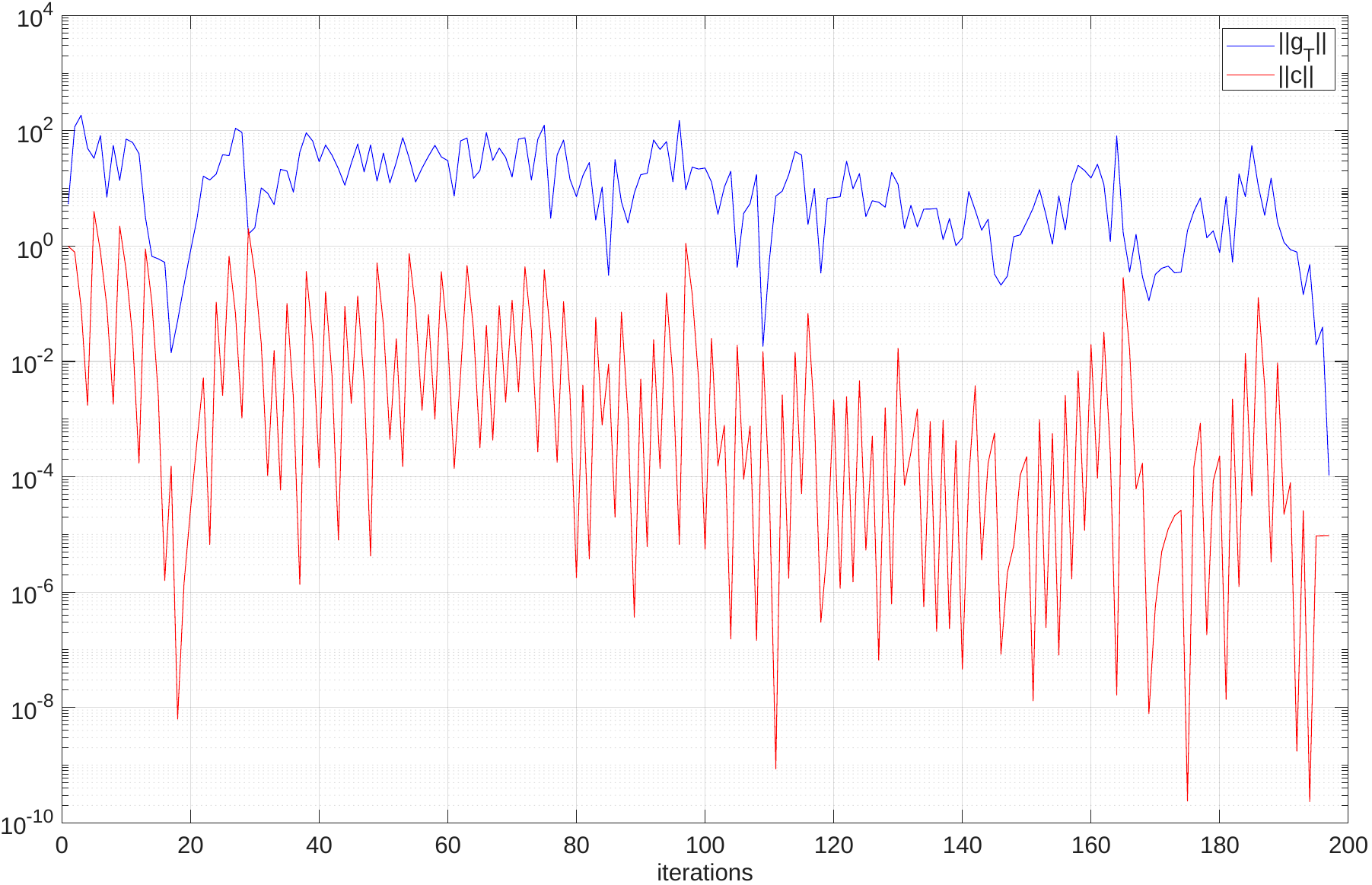}
\caption{\label{fig:noisy}Evolution of the two components of the
  optimality measure for a run of the noisy BT1 problem, with 5\% (top left),
  15\% (top right), 25\% (bottom left) and 50\% (bottom right)
  relative Gaussian noise in the objective function's gradient.  Despite severe random oscillations,
  the optimality measure shows the downward trend predicted by the
  deterministic theory.}
\end{figure}

To illustrate the resilience of \al{ADSWITCH} when noise is
present in the objective function's gradient, we also ran the algorithm on all instances of our test
problems while adding relative Gaussian noise (of zero mean and unit
variance) to the gradient. Each test problem was run ten times
independently, and for four levels (5\%, 15\%, 25\% and 50\%) of
noise in the objective function's gradient. The termination $\epsilon$ was set to $10^{-3}$ in
\req{convg}--\req{accept}.  The table in Figure~\ref{tab:reliab}
reports, for each of these levels, the number of problems for which
all ten runs failed (middle column) and the number of problems for
which all ten runs converged (rightmost column).  This table shows
that the reliability of the algorithm is remarkably insensitive to the
noise level. It is indeed noteworthy that around two thirds of the
considered test problems could be solved with top reliability even if
their gradients are perturbed by 50\% relative noise, which means that
barely one significant digit is correct.  This stable behaviour is
also visible in Figure~\ref{fig:noisy} showing the performance of the
algorithm on the {\tt BT1} problem for increasing noise levels.
  
\numsection{Conclusions and perspectives}\label{sec:concl}

We have proposed \al{ADSWITCH}, a very simple first-order algorithm for solving
nonlinear optimization problems with nonlinear equality
constraints. \al{ADSWITCH} adaptively selects steps in the plane
tangent to the constraints or steps that reduce infeasibility. It does
so without using a merit function or a filter to enforce convergence,
but instead relies on the simple switching condition \req{switch}. The
tangent steps are based on the AdaGrad method for unconstrained
minimization. As with AdaGrad, the objective function is
never evaluated.  We have analyzed its worst-case evaluation
complexity and obtain a global convergence rate that matches the best
known rates for unconstrained problems. Numerical experiments have
been presented indicating that the performance of \al{ADSWITCH} is
comparable to that of first-order methods for unconstrained problems.

At this stage, several theoretical questions remain open for further
research. They include a stochastic convergence analysis, 
the cases of rank-deficient Jacobians 
and/or unbounded gradients, the use of alternative minimization
methods in the tangential step (such as \al{ADAM} \cite{KingBa15},
\al{ASTR1} \cite{GratJeraToin22b}, \al{OFFAR}
\cite{GratJeraToin24b} or more standard unconstrained minimizing
techniques using objective function values).
On the practical side, further numerical experiments are clearly
necessary for assessing the true potential of the new method.  The
handling of inequality constraints is also of obvious interest.

{\footnotesize

\section*{\footnotesize Acknowledgement}

Philippe Toint is grateful for the continued and
friendly support of the Algo team at Toulouse IRIT. The authors also
thank Benedetta Morini and Stefania Bellavia for pointing out a number
of typos. Finally, the authors are indebted to an anonymous referee
for his/her very careful reading, which led to significant improvements.

}

\appendix

\appnumsection{Appendix; Details of the Numerical Results}

The following tables give details of the numerical results
discussed in Section~\ref{sec:numerics}.
\begin{itemize}
\item Table~\ref{no-noise} reports the results of the noiseless
  runs. Columns 2 and 3 indicate the number of variables and the
  number of constraints. Column ``\#its'' gives the number of iterations and 
  column ``exitc'' the termination condition (``convg'' = convergence,
  ``infeas'' = infeasible critical point, ``maxit'' = maximum number
  of iterations reached).
\item Tables~\ref{noise5}, \ref{noise15}, \ref{noise25} and
  \ref{noise50} report the average statistics of 10 independent runs
  for each problem, for increasing relative Gaussian noise levels (in
  the gradient of the objective function).
  The last column of each table gives the number of successful runs
  (among 10).
\end{itemize}

\begin{table}{\tiny
\begin{center}
\caption{\label{no-noise}Results of running \al{ADSWITCH} in the absence of noise}
\begin{tabular}{|l|r|r|r|r|r|r|r|}
\hline
Problem     &  $n$ & $m$ &    $f(x)$~~~~~     &  $\|g_T(x)\|$ & $\|c(x)\|$ & \#its & exitc \\
\hline
BT1        &    2 &    1 &  -9.999918e-01 & 9.64e-07 & 8.25e-08 &    141 & convg \\
BT2        &    3 &    1 &  +3.256821e-02 & 9.34e-07 & 6.29e-07 &    186 & convg \\
BT3        &    5 &    3 &  +4.093023e+00 & 8.58e-07 & 1.76e-08 &     70 & convg \\
BT4        &    3 &    2 &  -4.551055e+01 & 4.49e-07 & 1.16e-07 &     16 & convg \\
BT5        &    3 &    2 &  +9.617152e+02 & 7.00e-07 & 1.86e-07 &     46 & convg \\
BT6        &    5 &    2 &  +2.770448e-01 & 8.75e-07 & 1.55e-08 &    107 & convg \\
BT7        &    5 &    3 &      ---       &   ---    &   ---    &  10572 & infeas \\
BT8        &    5 &    2 &      ---       &   ---    &   ---    &     21 & infeas \\
BT9        &    4 &    2 &  -1.000000e+00 & 9.95e-07 & 4.92e-07 &   2481 & convg \\
BT10       &    2 &    2 &  -1.000000e+00 & 0.00e+00 & 5.85e-09 &      6 & convg \\
BT11       &    5 &    3 &  +8.248909e-01 & 9.55e-07 & 7.76e-07 &    351 & convg \\
BT12       &    5 &    3 &  +6.188119e+00 & 9.46e-07 & 4.53e-07 &    240 & convg \\
BYRDSPHR   &    3 &    2 &  -4.683300e+00 & 8.94e-07 & 3.33e-08 &     50 & convg \\
DIXCHLNG   &   10 &    5 &  +2.471898e+03 & 9.98e-07 & 2.67e-07 &   2722 & convg \\
EIGENA2    &  110 &   55 &  +2.448710e-13 & 9.90e-07 & 2.06e-09 &    696 & convg \\
EIGENACO   &  110 &   55 &  +2.436158e-13 & 9.87e-07 & 9.41e-16 &    389 & convg \\
EIGENB2    &  110 &   55 &  +1.800000e+01 & 6.82e-07 & 1.78e-15 &     18 & convg \\
EIGENBCO   &  110 &   55 &  +9.000000e+00 & 6.82e-07 & 1.66e-15 &     18 & convg \\
ELEC       &   75 &   25 &  +2.438128e+02 & 1.00e-06 & 4.79e-07 &  65380 & convg \\
GENHS28    &   10 &    8 &  +9.271737e-01 & 9.27e-07 & 1.36e-11 &     85 & convg \\
HS100LNP   &    7 &    2 &  +6.806301e+02 & 8.57e-07 & 2.47e-07 &    113 & convg \\
HS6        &    2 &    1 &  +9.238528e-13 & 8.60e-07 & 8.53e-07 &     55 & convg \\
HS7        &    2 &    1 &  -1.732051e+00 & 9.26e-07 & 1.96e-07 &    284 & convg \\
HS8        &    2 &    2 &  -1.000000e+00 & 0.00e+00 & 1.32e-08 &      4 & convg \\
HS9        &    2 &    1 &  -5.000000e-01 & 8.87e-07 & 7.11e-15 &     22 & convg \\
HS26       &    3 &    1 &  +1.303916e-09 & 1.00e-06 & 5.10e-07 &  18098 & convg \\
HS27       &    3 &    1 &  +3.999998e-02 & 9.79e-07 & 5.23e-07 &    271 & convg \\
HS28       &    3 &    1 &  +9.687529e-13 & 9.02e-07 & 1.33e-15 &    137 & convg \\
HS39       &    4 &    2 &  -1.000000e+00 & 9.95e-07 & 4.92e-07 &   2481 & convg \\
HS40       &    4 &    3 &  -2.500002e-01 & 9.74e-07 & 3.24e-07 &    600 & convg \\
HS42       &    4 &    2 &  +1.385786e+01 & 9.85e-07 & 8.75e-07 &     92 & convg \\
HS46       &    5 &    2 &  +8.808661e-09 & 1.00e-06 & 2.35e-07 &  13907 & convg \\
HS47       &    5 &    3 &  +4.650485e-10 & 1.00e-06 & 6.88e-07 &  26784 & convg \\
HS48       &    5 &    2 &  +2.912966e-13 & 9.31e-07 & 2.81e-15 &    177 & convg \\
HS50       &    5 &    3 &  +2.888537e-13 & 9.99e-07 & 6.63e-14 &  11166 & convg \\
HS51       &    5 &    3 &  +8.485715e-14 & 5.67e-07 & 1.46e-15 &     19 & convg \\
HS52       &    5 &    3 &  +5.326648e+00 & 9.87e-07 & 1.76e-09 &    150 & convg \\
HS61       &    3 &    2 &      ---       &   ---    &   ---    &      2 & infeas \\
HS77       &    5 &    2 &  +2.415051e-01 & 9.27e-07 & 1.18e-07 &    103 & convg \\
HS78       &    5 &    3 &  -2.919700e+00 & 9.31e-07 & 1.97e-08 &     95 & convg \\
HS79       &    5 &    3 &  +7.877683e-02 & 9.71e-07 & 2.29e-07 &    357 & convg \\
LUKVLE1    &   20 &   18 &  +6.232459e+00 & 9.11e-07 & 1.20e-08 &     98 & convg \\
LUKVLE2    &   20 &   13 &  +4.372199e+02 & 9.93e-07 & 3.19e-08 &    604 & convg \\
LUKVLE3    &   20 &    2 &  +2.758657e+01 & 9.94e-07 & 8.67e-07 &   2319 & convg \\
LUKVLE4    &   20 &    9 &  +1.060660e+02 & 1.36e+02 & 3.15e-07 & 100000 & maxit \\
LUKVLE6    &   21 &   10 &  +1.120010e+03 & 2.10e-01 & 7.64e-07 & 100000 & maxit \\
LUKVLE7    &   20 &    4 &  -5.352208e+00 & 8.90e-07 & 3.91e-08 &    111 & convg \\
LUKVLE8    &   20 &   18 &  +2.099621e+03 & 8.65e-07 & 3.18e-09 &     55 & convg \\
LUKVLE9    &   20 &    6 &  +6.560168e+00 & 1.51e+01 & 8.57e-07 & 100000 & maxit \\
LUKVLE10   &   20 &   18 &  +6.639332e+00 & 7.23e-07 & 3.27e-07 &     40 & convg \\
LUKVLE11   &   18 &   10 &  +5.359296e-07 & 3.53e-05 & 1.69e-09 & 100000 & maxit \\
LUKVLE12   &   17 &   12 &  +2.286872e+02 & 4.44e+01 & 2.23e+00 & 100000 & maxit \\
LUKVLE13   &   18 &   10 &  +5.483076e+01 & 4.24e-06 & 1.20e-08 & 100000 & maxit \\
LUKVLE14   &   18 &   10 &  +4.180657e+04 & 7.86e-01 & 7.02e-07 & 100000 & maxit \\
LUKVLE15   &   17 &   12 &  +7.527144e+01 & 1.14e+01 & 9.10e-07 & 100000 & maxit \\
LUKVLE16   &   17 &   12 &  +8.040898e+01 & 2.71e+01 & 3.79e+00 & 100000 & maxit \\
LUKVLE17   &   17 &   12 &  +1.190180e+02 & 5.48e+01 & 8.37e+00 & 100000 & maxit \\
LUKVLE18   &   17 &   12 &  +2.627594e+01 & 1.71e+01 & 8.37e+00 & 100000 & maxit \\
LUKVLI4    &   20 &    9 &  +1.060660e+02 & 1.36e+02 & 3.15e-07 & 100000 & maxit \\
MARATOS    &    2 &    1 &  -1.000000e+00 & 9.27e-07 & 2.26e-07 &    168 & convg \\
MWRIGHT    &    5 &    3 &  +2.497881e+01 & 9.23e-07 & 1.04e-08 &    114 & convg \\
ORTHRDM2   &    9 &    3 &  +1.087113e-13 & 6.51e-07 & 7.14e-07 &     31 & convg \\
ORTHRDS2   &    9 &    3 &  +5.921433e-14 & 4.81e-07 & 2.63e-08 &     33 & convg \\
ORTHREGA   &  133 &   64 &  +3.503002e+02 & 9.82e-07 & 1.66e-08 &    705 & convg \\
ORTHREGB   &   27 &    6 &  +3.057742e-14 & 3.50e-07 & 1.02e-08 &     15 & convg \\
ORTHREGC   &   15 &    5 &  +1.125393e-12 & 9.24e-07 & 7.22e-08 &     50 & convg \\
ORTHREGD   &   43 &   20 &  +2.179046e+02 & 1.74e+01 & 1.09e-03 & 100000 & maxit \\
ORTHRGDM   &   43 &   20 &  +5.067350e+01 & 8.20e+00 & 1.97e-03 & 100000 & maxit \\
ORTHRGDS   &   43 &   20 &  +6.210193e+00 & 9.97e-07 & 8.69e-07 &    474 & convg \\
S316m322   &    2 &    1 &      ---       &   ---    &   ---    &      0 & infeas \\
SPINOP     &   11 &    9 &  +8.435060e-02 & 7.91e-04 & 2.20e-06 & 100000 & maxit \\
\hline
\end{tabular}
\end{center}
}
\end{table}

\begin{table}{\tiny
\begin{center}
\caption{\label{noise5}Results of running \al{ADSWITCH} with 5\% relative Gaussian noise}
\begin{tabular}{|l|r|r|r|r|r|r|r|}
\hline
Problem     &  $n$ & $m$ &  avr $f(x)$ & avr $\|g_T(x)\|$ & avr $\|c(x)\|$ & avr \#its & \#success \\
\hline
BT1        &    3 &    1 &  -5.681296e-01 & 7.73e-04 & 3.20e-04 & 8.92e+01 & 10 \\
BT2        &    5 &    3 &  +3.257271e-02 & 1.51e-03 & 3.49e-04 & 1.18e+02 & 10 \\
BT3        &    3 &    2 &  +4.089376e+00 & 5.05e-02 & 6.34e-04 & 2.50e+04 & 10 \\
BT4        &    3 &    2 &  -4.551215e+01 & 1.87e-02 & 5.03e-04 & 8.88e+01 & 10 \\
BT5        &    5 &    2 &  +9.617275e+02 & 1.60e-01 & 9.36e-05 & 1.10e+03 & 10 \\
BT6        &    5 &    3 &  +2.770621e-01 & 9.19e-03 & 2.69e-04 & 3.50e+03 & 10 \\
BT7        &    5 &    2 &      ---       &   ---    &   ---    & 2.51e+03 & 10 \\
BT8        &    4 &    2 &  +1.000238e+00 & 6.34e-04 & 3.41e-04 & 1.47e+01 & 10 \\
BT9        &    2 &    2 &  -1.000535e+00 & 9.87e-04 & 5.37e-04 & 1.15e+03 & 10 \\
BT10       &    5 &    3 &  -1.000038e+00 & 0.00e+00 & 5.91e-05 & 5.00e+00 & 10 \\
BT11       &    5 &    3 &  +8.246795e-01 & 4.64e-02 & 4.07e-04 & 1.29e+04 & 10 \\
BT12       &    3 &    2 &  +6.188179e+00 & 6.06e-03 & 6.49e-04 & 1.29e+02 & 10 \\
BYRDSPHR   &   10 &    5 &  -4.683084e+00 & 1.21e-02 & 1.63e-05 & 1.85e+02 & 10 \\
DIXCHLNG   &  110 &   55 &  +2.471678e+03 & 1.25e+01 & 3.47e-04 &    ---   &  9 \\
EIGENA2    &  110 &   55 &  +2.507664e-07 & 1.00e-03 & 4.03e-04 & 3.98e+02 & 10 \\
EIGENACO   &  110 &   55 &  +2.410396e-07 & 9.81e-04 & 9.35e-16 & 2.31e+02 & 10 \\
EIGENB2    &  110 &   55 &  +3.530739e-06 & 1.00e-03 & 6.16e-04 & 1.10e+03 & 10 \\
EIGENBCO   &   75 &   25 &  +5.356264e-06 & 9.58e-04 & 5.29e-04 & 1.07e+03 & 10 \\
ELEC       &   10 &    8 &  +2.438048e+02 & 1.86e-01 & 4.33e-04 & 1.00e+05 & 10 \\
GENHS28    &    7 &    2 &  +9.274072e-01 & 2.82e-02 & 6.75e-06 & 3.01e+03 & 10 \\
HS100LNP   &    2 &    1 &  +6.806369e+02 & 5.73e-01 & 1.91e-04 &    ---   &  1 \\
HS6        &    2 &    1 &  +9.637182e-07 & 8.72e-04 & 2.76e-04 & 2.96e+01 & 10 \\
HS7        &    2 &    2 &  -1.732173e+00 & 9.72e-04 & 4.26e-04 & 1.65e+02 & 10 \\
HS8        &    2 &    1 &  -1.000000e+00 & 0.00e+00 & 1.32e-08 & 4.00e+00 & 10 \\
HS9        &    3 &    1 &  -4.998315e-01 & 3.11e-03 & 5.86e-15 & 1.41e+01 & 10 \\
HS26       &    3 &    1 &  +1.374024e-05 & 1.01e-03 & 7.19e-04 & 2.19e+02 & 10 \\
HS27       &    3 &    1 &  +3.999576e-02 & 9.84e-04 & 2.58e-04 & 6.39e+01 & 10 \\
HS28       &    4 &    2 &  +1.022606e-06 & 9.25e-04 & 6.55e-16 & 7.09e+01 & 10 \\
HS39       &    4 &    3 &  -1.000759e+00 & 9.86e-04 & 7.61e-04 & 1.08e+03 & 10 \\
HS40       &    4 &    2 &  -2.501555e-01 & 1.42e-02 & 4.17e-04 & 2.71e+02 & 10 \\
HS42       &    5 &    2 &  +1.386050e+01 & 1.79e-01 & 2.38e-04 & 9.08e+02 & 10 \\
HS46       &    5 &    3 &  +4.853021e-06 & 9.96e-04 & 5.22e-04 & 1.31e+02 & 10 \\
HS47       &    5 &    2 &  +1.474970e-05 & 1.03e-03 & 4.84e-04 & 7.36e+02 & 10 \\
HS48       &    5 &    3 &  +3.056976e-07 & 9.53e-04 & 1.88e-15 & 9.29e+01 & 10 \\
HS50       &    5 &    3 &  +3.332153e-07 & 1.07e-03 & 1.07e-13 & 6.29e+03 & 10 \\
HS51       &    5 &    3 &  +1.176429e-07 & 6.55e-04 & 1.69e-15 & 1.03e+01 & 10 \\
HS52       &    3 &    2 &  +5.326494e+00 & 1.14e-01 & 6.34e-05 & 3.00e+04 & 10 \\
HS61       &    5 &    2 &      ---       &   ---    &   ---    & 1.00e+00 & 10 \\
HS77       &    5 &    3 &  +2.415176e-01 & 1.14e-02 & 3.55e-04 & 2.04e+03 & 10 \\
HS78       &    5 &    3 &  -2.918759e+00 & 6.58e-02 & 3.58e-04 &    ---   &  6 \\
HS79       &   20 &   18 &  +7.879160e-02 & 2.26e-03 & 3.19e-04 & 1.24e+02 & 10 \\
LUKVLE1    &   20 &   13 &  +6.234044e+00 & 2.71e-03 & 3.84e-04 & 1.13e+02 & 10 \\
LUKVLE2    &   20 &    2 &  +4.371168e+02 & 3.46e+00 & 5.67e-04 & 1.00e+05 & 10 \\
LUKVLE3    &   20 &    9 &  +2.757471e+01 & 1.64e-01 & 6.39e-04 & 1.00e+05 & 10 \\
LUKVLE4    &   21 &   10 &  +1.058513e+02 & 1.35e+02 & 5.50e-06 & 1.00e+05 & 10 \\
LUKVLE6    &   20 &    4 &  +1.119942e+03 & 4.69e-01 & 5.08e-04 & 1.00e+05 & 10 \\
LUKVLE7    &   20 &   18 &  -5.353248e+00 & 1.18e-01 & 3.97e-04 &    ---   &  0 \\
LUKVLE8    &   20 &    6 &      ---       &   ---    &   ---    & 1.13e+04 & 10 \\
LUKVLE9    &   20 &   18 &  +6.555645e+00 & 1.51e+01 & 1.66e-04 & 1.00e+05 & 10 \\
LUKVLE10   &   18 &   10 &  +6.639743e+00 & 1.44e-02 & 2.76e-04 & 4.84e+03 & 10 \\
LUKVLE11   &   17 &   12 &  +3.438763e-05 & 1.09e-03 & 4.80e-04 & 6.96e+03 & 10 \\
LUKVLE12   &   18 &   10 &  +2.286872e+02 & 4.44e+01 & 2.23e+00 & 1.00e+05 & 10 \\
LUKVLE13   &   18 &   10 &  +5.483466e+01 & 8.14e-02 & 4.72e-04 & 1.00e+05 & 10 \\
LUKVLE14   &   17 &   12 &  +4.180658e+04 & 8.46e-01 & 7.29e-04 & 1.00e+05 & 10 \\
LUKVLE15   &   17 &   12 &  +7.511104e+01 & 1.14e+01 & 1.81e-04 &    ---   &  0 \\
LUKVLE16   &   17 &   12 &  +8.040898e+01 & 2.71e+01 & 3.79e+00 &    ---   &  0 \\
LUKVLE17   &   17 &   12 &  +1.190180e+02 & 5.48e+01 & 8.37e+00 &    ---   &  0 \\
LUKVLE18   &   20 &    9 &  +2.627594e+01 & 1.71e+01 & 8.37e+00 &    ---   &  0 \\
LUKVLI4    &    2 &    1 &  +1.063358e+02 & 1.38e+02 & 5.60e-06 & 1.00e+05 & 10 \\
MARATOS    &    5 &    3 &  -1.000295e+00 & 9.53e-04 & 5.91e-04 & 9.98e+01 & 10 \\
MWRIGHT    &    9 &    3 &  +2.497543e+01 & 1.45e-01 & 4.88e-04 &    ---   &  9 \\
ORTHRDM2   &    9 &    3 &  +1.034714e-07 & 6.19e-04 & 3.51e-04 & 2.16e+01 & 10 \\
ORTHRDS2   &  133 &   64 &  +8.522014e-08 & 5.48e-04 & 3.05e-04 & 2.14e+01 & 10 \\
ORTHREGA   &   27 &    6 &  +3.502947e+02 & 6.49e-02 & 7.05e-04 & 1.00e+05 & 10 \\
ORTHREGB   &   15 &    5 &  +5.311601e-08 & 3.99e-04 & 3.80e-05 & 1.06e+01 & 10 \\
ORTHREGC   &   43 &   20 &  +6.775424e-07 & 8.16e-04 & 5.67e-04 & 2.69e+01 & 10 \\
ORTHREGD   &   43 &   20 &      ---       &   ---    &   ---    & 2.13e+02 & 10 \\
ORTHRGDM   &   43 &   20 &      ---       &   ---    &   ---    & 1.86e+03 & 10 \\
ORTHRGDS   &    2 &    1 &  +6.210324e+00 & 2.21e-02 & 5.60e-04 &    ---   &  0 \\
S316m322   &   11 &    9 &      ---       &   ---    &   ---    & 0.00e+00 & 10 \\
SPINOP     &   11 &    9 &  +9.941025e-02 & 1.10e-03 & 5.74e-04 & 6.14e+04 & 10 \\
\hline
\end{tabular}
\end{center}
}
\end{table}
\begin{table}{\tiny
\begin{center}
\caption{\label{noise15}Results of running \al{ADSWITCH} with 15\% relative Gaussian noise}
\begin{tabular}{|l|r|r|r|r|r|r|r|}
\hline
Problem     &  $n$ & $m$ &  avr $f(x)$ & avr $\|g_T(x)\|$ & avr $\|c(x)\|$ & avr \#its & \#success \\
\hline
BT1        &    3 &    1 &  -1.720654e-01 & 2.87e-04 & 2.80e-04 & 2.19e+02 & 10 \\
BT2        &    5 &    3 &  +3.257574e-02 & 2.86e-03 & 5.14e-04 & 1.45e+02 & 10 \\
BT3        &    3 &    2 &  +4.089932e+00 & 6.74e-02 & 6.34e-04 & 7.06e+04 & 10 \\
BT4        &    3 &    2 &  -4.551083e+01 & 5.03e-02 & 1.67e-04 & 3.36e+02 & 10 \\
BT5        &    5 &    2 &  +9.617453e+02 & 2.74e-01 & 4.39e-05 & 4.57e+03 & 10 \\
BT6        &    5 &    3 &  +2.771418e-01 & 1.71e-02 & 3.62e-04 &    ---   &  7 \\
BT7        &    5 &    2 &      ---       &   ---    &   ---    & 1.97e+03 & 10 \\
BT8        &    4 &    2 &  +1.000238e+00 & 3.77e-04 & 3.41e-04 & 1.38e+01 & 10 \\
BT9        &    2 &    2 &  -1.000468e+00 & 7.93e-04 & 4.70e-04 & 6.71e+02 & 10 \\
BT10       &    5 &    3 &  -1.000038e+00 & 0.00e+00 & 5.91e-05 & 5.00e+00 & 10 \\
BT11       &    5 &    3 &  +8.248886e-01 & 6.47e-02 & 4.69e-04 &    ---   &  7 \\
BT12       &    3 &    2 &  +6.188227e+00 & 1.06e-02 & 4.39e-04 & 1.90e+02 & 10 \\
BYRDSPHR   &   10 &    5 &  -4.680516e+00 & 3.82e-02 & 1.76e-04 & 4.58e+02 & 10 \\
DIXCHLNG   &  110 &   55 &  +2.471747e+03 & 3.29e+01 & 5.13e-04 &    ---   &  8 \\
EIGENA2    &  110 &   55 &  +2.864156e-07 & 1.07e-03 & 5.67e-04 & 3.97e+02 & 10 \\
EIGENACO   &  110 &   55 &  +2.500993e-07 & 9.99e-04 & 9.37e-16 & 2.35e+02 & 10 \\
EIGENB2    &  110 &   55 &  +2.989956e-06 & 1.02e-03 & 2.20e-04 & 1.05e+03 & 10 \\
EIGENBCO   &   75 &   25 &  +1.674676e-06 & 8.43e-04 & 4.82e-04 & 1.35e+03 & 10 \\
ELEC       &   10 &    8 &  +2.438072e+02 & 4.00e-01 & 5.36e-04 &    ---   &  6 \\
GENHS28    &    7 &    2 &  +9.275303e-01 & 3.38e-02 & 6.75e-06 & 3.59e+04 & 10 \\
HS100LNP   &    2 &    1 &  +6.806511e+02 & 1.04e+00 & 3.81e-04 &    ---   &  0 \\
HS6        &    2 &    1 &  +8.822774e-07 & 8.34e-04 & 5.77e-04 & 2.94e+01 & 10 \\
HS7        &    2 &    2 &  -1.732173e+00 & 9.28e-04 & 4.25e-04 & 1.45e+02 & 10 \\
HS8        &    2 &    1 &  -1.000000e+00 & 0.00e+00 & 1.32e-08 & 4.00e+00 & 10 \\
HS9        &    3 &    1 &  -4.993365e-01 & 7.17e-03 & 6.76e-15 & 2.20e+01 & 10 \\
HS26       &    3 &    1 &  +1.632328e-05 & 1.15e-03 & 4.97e-04 & 2.15e+02 & 10 \\
HS27       &    3 &    1 &  +3.999706e-02 & 9.14e-04 & 2.16e-04 & 4.66e+01 & 10 \\
HS28       &    4 &    2 &  +1.738691e-06 & 1.20e-03 & 7.54e-16 & 6.87e+01 & 10 \\
HS39       &    4 &    3 &  -1.000411e+00 & 1.07e-03 & 4.13e-04 & 5.98e+02 & 10 \\
HS40       &    4 &    2 &  -2.493358e-01 & 3.90e-02 & 4.84e-06 & 4.18e+02 & 10 \\
HS42       &    5 &    2 &  +1.386567e+01 & 2.88e-01 & 1.42e-04 & 1.38e+03 & 10 \\
HS46       &    5 &    3 &  +7.154893e-06 & 1.06e-03 & 3.65e-04 & 1.24e+02 & 10 \\
HS47       &    5 &    2 &  +1.763663e-05 & 1.17e-03 & 5.01e-04 & 6.97e+02 & 10 \\
HS48       &    5 &    3 &  +3.810065e-07 & 1.05e-03 & 1.86e-15 & 9.47e+01 & 10 \\
HS50       &    5 &    3 &  +4.796120e-07 & 1.29e-03 & 9.03e-14 & 6.20e+03 & 10 \\
HS51       &    5 &    3 &  +1.153904e-07 & 6.35e-04 & 1.57e-15 & 9.90e+00 & 10 \\
HS52       &    3 &    2 &  +5.326745e+00 & 1.30e-01 & 6.34e-05 &    ---   &  5 \\
HS61       &    5 &    2 &      ---       &   ---    &   ---    & 1.00e+00 & 10 \\
HS77       &    5 &    3 &  +2.415599e-01 & 1.79e-02 & 2.39e-04 & 5.11e+04 & 10 \\
HS78       &    5 &    3 &  -2.918552e+00 & 6.92e-02 & 4.09e-04 &    ---   &  4 \\
HS79       &   20 &   18 &  +7.879634e-02 & 2.28e-03 & 5.15e-04 & 1.54e+02 & 10 \\
LUKVLE1    &   20 &   13 &  +6.233718e+00 & 7.71e-03 & 3.12e-04 & 1.20e+02 & 10 \\
LUKVLE2    &   20 &    2 &  +4.371783e+02 & 7.13e+00 & 4.04e-04 &    ---   &  8 \\
LUKVLE3    &   20 &    9 &  +2.757958e+01 & 5.72e-01 & 4.40e-04 &    ---   &  9 \\
LUKVLE4    &   21 &   10 &  +1.063795e+02 & 1.41e+02 & 8.44e-06 & 1.00e+05 & 10 \\
LUKVLE6    &   20 &    4 &  +1.119853e+03 & 1.16e+00 & 5.46e-04 & 1.00e+05 & 10 \\
LUKVLE7    &   20 &   18 &  -5.350163e+00 & 4.22e-01 & 2.45e-04 &    ---   &  0 \\
LUKVLE8    &   20 &    6 &      ---       &   ---    &   ---    & 1.12e+03 & 10 \\
LUKVLE9    &   20 &   18 &  +6.707674e+00 & 1.60e+01 & 2.84e-04 & 1.00e+05 & 10 \\
LUKVLE10   &   18 &   10 &  +6.473732e+00 & 2.91e-02 & 3.32e-04 & 1.19e+04 & 10 \\
LUKVLE11   &   17 &   12 &  +2.849669e-01 & 5.36e-03 & 5.57e-04 &    ---   &  9 \\
LUKVLE12   &   18 &   10 &  +2.286872e+02 & 4.44e+01 & 2.23e+00 & 1.00e+05 & 10 \\
LUKVLE13   &   18 &   10 &  +5.483592e+01 & 1.08e-01 & 5.89e-04 & 1.00e+05 & 10 \\
LUKVLE14   &   17 &   12 &  +4.180658e+04 & 1.02e+00 & 3.89e-04 & 1.00e+05 & 10 \\
LUKVLE15   &   17 &   12 &  +7.514427e+01 & 1.15e+01 & 1.81e-04 &    ---   &  0 \\
LUKVLE16   &   17 &   12 &  +8.040898e+01 & 2.71e+01 & 3.79e+00 &    ---   &  0 \\
LUKVLE17   &   17 &   12 &  +1.190180e+02 & 5.48e+01 & 8.37e+00 &    ---   &  0 \\
LUKVLE18   &   20 &    9 &  +2.627594e+01 & 1.71e+01 & 8.37e+00 &    ---   &  0 \\
LUKVLI4    &    2 &    1 &  +1.094837e+02 & 1.56e+02 & 6.13e-06 & 1.00e+05 & 10 \\
MARATOS    &    5 &    3 &  -1.000232e+00 & 7.86e-04 & 4.65e-04 & 9.35e+01 & 10 \\
MWRIGHT    &    9 &    3 &  +2.497649e+01 & 2.93e-01 & 5.48e-04 &    ---   &  7 \\
ORTHRDM2   &    9 &    3 &  +9.737407e-08 & 5.56e-04 & 4.05e-04 & 2.27e+01 & 10 \\
ORTHRDS2   &  133 &   64 &  +1.610065e-07 & 7.37e-04 & 4.02e-04 & 2.24e+01 & 10 \\
ORTHREGA   &   27 &    6 &  +3.502994e+02 & 1.19e-01 & 4.88e-04 &    ---   &  5 \\
ORTHREGB   &   15 &    5 &  +7.691602e-08 & 4.91e-04 & 8.40e-05 & 1.09e+01 & 10 \\
ORTHREGC   &   43 &   20 &  +9.146250e-07 & 8.78e-04 & 3.23e-04 & 2.67e+01 & 10 \\
ORTHREGD   &   43 &   20 &      ---       &   ---    &   ---    & 1.41e+02 & 10 \\
ORTHRGDM   &   43 &   20 &      ---       &   ---    &   ---    & 1.66e+03 & 10 \\
ORTHRGDS   &    2 &    1 &  +6.210852e+00 & 5.11e-02 & 5.33e-04 &    ---   &  0 \\
S316m322   &   11 &    9 &      ---       &   ---    &   ---    & 0.00e+00 & 10 \\
SPINOP     &   11 &    9 &  +1.134528e-01 & 1.43e-03 & 6.21e-04 & 4.12e+04 & 10 \\
\hline
\end{tabular}
\end{center}
}
\end{table}
\begin{table}{\tiny
\begin{center}
\caption{\label{noise25}Results of running \al{ADSWITCH} with 25\% relative Gaussian noise}
\begin{tabular}{|l|r|r|r|r|r|r|r|}
\hline
Problem     &  $n$ & $m$ &  avr $f(x)$ & avr $\|g_T(x)\|$ & avr $\|c(x)\|$ & avr \#its & \#success \\
\hline
BT1        &    3 &    1 &  -1.628888e-01 & 1.65e-05 & 3.72e-04 & 1.86e+02 & 10 \\
BT2        &    5 &    3 &  +3.257851e-02 & 5.28e-03 & 3.90e-04 & 2.48e+02 & 10 \\
BT3        &    3 &    2 &  +4.089948e+00 & 6.60e-02 & 6.34e-04 & 9.29e+04 & 10 \\
BT4        &    3 &    2 &  -4.551030e+01 & 1.04e-01 & 2.53e-04 & 5.51e+02 & 10 \\
BT5        &    5 &    2 &  +9.617367e+02 & 2.09e-01 & 1.93e-04 & 6.10e+03 & 10 \\
BT6        &    5 &    3 &  +2.771881e-01 & 2.53e-02 & 6.25e-04 &    ---   &  1 \\
BT7        &    5 &    2 &      ---       &   ---    &   ---    & 1.51e+03 & 10 \\
BT8        &    4 &    2 &  +1.000238e+00 & 2.43e-04 & 3.41e-04 & 1.42e+01 & 10 \\
BT9        &    2 &    2 &  -1.000286e+00 & 6.60e-04 & 2.91e-04 & 4.16e+02 & 10 \\
BT10       &    5 &    3 &  -1.000038e+00 & 0.00e+00 & 5.91e-05 & 5.00e+00 & 10 \\
BT11       &    5 &    3 &  +8.259372e-01 & 1.09e-01 & 4.43e-04 &    ---   &  6 \\
BT12       &    3 &    2 &  +6.188639e+00 & 2.01e-02 & 3.28e-04 & 2.05e+02 & 10 \\
BYRDSPHR   &   10 &    5 &  -4.672953e+00 & 8.00e-02 & 9.61e-05 & 5.88e+02 & 10 \\
DIXCHLNG   &  110 &   55 &  +2.471810e+03 & 3.74e+01 & 4.88e-04 &    ---   &  6 \\
EIGENA2    &  110 &   55 &  +2.814117e-07 & 1.06e-03 & 3.87e-04 & 4.03e+02 & 10 \\
EIGENACO   &  110 &   55 &  +2.651002e-07 & 1.03e-03 & 9.35e-16 & 2.43e+02 & 10 \\
EIGENB2    &  110 &   55 &  +2.653621e-06 & 1.02e-03 & 5.90e-04 & 1.27e+03 & 10 \\
EIGENBCO   &   75 &   25 &  +6.332410e-07 & 8.51e-04 & 4.14e-04 & 1.56e+03 & 10 \\
ELEC       &   10 &    8 &  +2.438104e+02 & 5.26e-01 & 5.69e-04 &    ---   &  6 \\
GENHS28    &    7 &    2 &  +9.281926e-01 & 5.42e-02 & 6.75e-06 & 2.64e+04 & 10 \\
HS100LNP   &    2 &    1 &  +6.806694e+02 & 1.50e+00 & 6.92e-04 &    ---   &  0 \\
HS6        &    2 &    1 &  +9.704023e-07 & 8.32e-04 & 4.74e-04 & 3.12e+01 & 10 \\
HS7        &    2 &    2 &  -1.732196e+00 & 7.94e-04 & 5.04e-04 & 1.51e+02 & 10 \\
HS8        &    2 &    1 &  -1.000000e+00 & 0.00e+00 & 1.32e-08 & 4.00e+00 & 10 \\
HS9        &    3 &    1 &  -4.981924e-01 & 1.09e-02 & 2.06e-14 & 7.99e+01 & 10 \\
HS26       &    3 &    1 &  +1.680639e-05 & 1.18e-03 & 5.93e-04 & 2.00e+02 & 10 \\
HS27       &    3 &    1 &  +4.000247e-02 & 9.90e-04 & 1.26e-04 & 4.02e+01 & 10 \\
HS28       &    4 &    2 &  +2.764455e-06 & 1.42e-03 & 8.77e-16 & 7.06e+01 & 10 \\
HS39       &    4 &    3 &  -1.000434e+00 & 7.07e-04 & 4.45e-04 & 3.96e+02 & 10 \\
HS40       &    4 &    2 &  -2.232808e-01 & 6.51e-02 & 1.23e-04 & 5.08e+02 & 10 \\
HS42       &    5 &    2 &  +1.387811e+01 & 4.89e-01 & 3.55e-04 & 3.79e+03 & 10 \\
HS46       &    5 &    3 &  +7.617645e-06 & 1.14e-03 & 6.62e-04 & 1.16e+02 & 10 \\
HS47       &    5 &    2 &  +2.249967e-05 & 1.37e-03 & 3.85e-04 & 6.16e+02 & 10 \\
HS48       &    5 &    3 &  +5.445390e-07 & 1.25e-03 & 1.39e-15 & 8.88e+01 & 10 \\
HS50       &    5 &    3 &  +7.538284e-07 & 1.61e-03 & 9.79e-14 & 6.08e+03 & 10 \\
HS51       &    5 &    3 &  +1.182718e-07 & 6.25e-04 & 1.67e-15 & 1.05e+01 & 10 \\
HS52       &    3 &    2 &  +5.327150e+00 & 1.37e-01 & 6.34e-05 &    ---   &  1 \\
HS61       &    5 &    2 &      ---       &   ---    &   ---    & 1.00e+00 & 10 \\
HS77       &    5 &    3 &  +2.416327e-01 & 2.47e-02 & 2.62e-04 &    ---   &  8 \\
HS78       &    5 &    3 &  -2.918216e+00 & 8.49e-02 & 3.36e-04 &    ---   &  1 \\
HS79       &   20 &   18 &  +7.879582e-02 & 3.00e-03 & 6.26e-04 & 2.09e+02 & 10 \\
LUKVLE1    &   20 &   13 &  +4.363520e+00 & 6.83e-03 & 3.41e-04 & 1.25e+02 & 10 \\
LUKVLE2    &   20 &    2 &  +4.371531e+02 & 5.88e+00 & 6.31e-04 &    ---   &  9 \\
LUKVLE3    &   20 &    9 &  +2.758052e+01 & 1.14e+00 & 6.23e-04 &    ---   &  7 \\
LUKVLE4    &   21 &   10 &  +1.146093e+02 & 1.73e+02 & 6.14e-06 & 1.00e+05 & 10 \\
LUKVLE6    &   20 &    4 &  +1.119798e+03 & 2.06e+00 & 4.77e-04 & 1.00e+05 & 10 \\
LUKVLE7    &   20 &   18 &  -5.349869e+00 & 4.53e-01 & 3.32e-04 &    ---   &  0 \\
LUKVLE8    &   20 &    6 &      ---       &   ---    &   ---    & 3.95e+02 & 10 \\
LUKVLE9    &   20 &   18 &  +6.467339e+00 & 1.45e+01 & 1.53e-04 & 1.00e+05 & 10 \\
LUKVLE10   &   18 &   10 &  +6.464724e+00 & 2.88e-02 & 5.32e-04 & 4.79e+04 & 10 \\
LUKVLE11   &   17 &   12 &  +6.422930e-05 & 1.84e-03 & 3.02e-04 & 5.17e+03 & 10 \\
LUKVLE12   &   18 &   10 &  +2.286872e+02 & 4.44e+01 & 2.23e+00 & 1.00e+05 & 10 \\
LUKVLE13   &   18 &   10 &  +5.501071e+01 & 1.89e-01 & 3.69e-04 & 1.00e+05 & 10 \\
LUKVLE14   &   17 &   12 &  +4.180664e+04 & 1.77e+00 & 4.53e-04 & 1.00e+05 & 10 \\
LUKVLE15   &   17 &   12 &  +7.384004e+01 & 1.14e+01 & 1.82e-04 &    ---   &  0 \\
LUKVLE16   &   17 &   12 &  +8.040898e+01 & 2.71e+01 & 3.79e+00 &    ---   &  0 \\
LUKVLE17   &   17 &   12 &  +1.190180e+02 & 5.48e+01 & 8.37e+00 &    ---   &  0 \\
LUKVLE18   &   20 &    9 &  +2.627594e+01 & 1.71e+01 & 8.37e+00 &    ---   &  0 \\
LUKVLI4    &    2 &    1 &  +1.094630e+02 & 1.55e+02 & 9.07e-06 & 1.00e+05 & 10 \\
MARATOS    &    5 &    3 &  -9.998270e-01 & 8.83e-03 & 2.92e-04 & 3.76e+01 & 10 \\
MWRIGHT    &    9 &    3 &  +2.497751e+01 & 2.90e-01 & 5.26e-04 &    ---   &  7 \\
ORTHRDM2   &    9 &    3 &  +1.408678e-07 & 6.56e-04 & 4.15e-04 & 2.25e+01 & 10 \\
ORTHRDS2   &  133 &   64 &  +1.867325e-07 & 7.67e-04 & 4.32e-04 & 2.25e+01 & 10 \\
ORTHREGA   &   27 &    6 &  +3.503036e+02 & 1.71e-01 & 4.19e-04 &    ---   &  3 \\
ORTHREGB   &   15 &    5 &  +6.680941e-08 & 4.66e-04 & 2.55e-04 & 1.21e+01 & 10 \\
ORTHREGC   &   43 &   20 &  +5.552793e-07 & 7.69e-04 & 3.77e-04 & 2.81e+01 & 10 \\
ORTHREGD   &   43 &   20 &      ---       &   ---    &   ---    & 1.44e+02 & 10 \\
ORTHRGDM   &   43 &   20 &      ---       &   ---    &   ---    & 1.60e+03 & 10 \\
ORTHRGDS   &    2 &    1 &  +6.211177e+00 & 6.17e-02 & 6.38e-04 &    ---   &  0 \\
S316m322   &   11 &    9 &      ---       &   ---    &   ---    & 0.00e+00 & 10 \\
SPINOP     &   11 &    9 &  +1.331340e-01 & 1.97e-03 & 5.60e-04 & 2.56e+04 & 10 \\
\hline
\end{tabular}
\end{center}
}
\end{table}
\begin{table}{\tiny
\begin{center}
\caption{\label{noise50}Results of running \al{ADSWITCH} with 50\% relative Gaussian noise}
\begin{tabular}{|l|r|r|r|r|r|r|r|}
\hline
Problem     &  $n$ & $m$ &  avr $f(x)$ & avr $\|g_T(x)\|$ & avr $\|c(x)\|$ & avr \#its & \#success \\
\hline
BT1        &    3 &    1 &  -8.595453e-01 & 9.99e-02 & 4.59e-04 & 2.67e+02 & 10 \\
BT2        &    5 &    3 &  +3.258385e-02 & 7.70e-03 & 4.76e-04 & 7.86e+02 & 10 \\
BT3        &    3 &    2 &  +4.090844e+00 & 1.06e-01 & 6.34e-04 &    ---   &  9 \\
BT4        &    3 &    2 &  -4.550471e+01 & 2.77e-01 & 2.28e-04 & 1.19e+03 & 10 \\
BT5        &    5 &    2 &  +9.618963e+02 & 6.57e-01 & 2.61e-04 & 5.57e+03 & 10 \\
BT6        &    5 &    3 &  +2.773363e-01 & 2.99e-02 & 5.98e-04 &    ---   &  0 \\
BT7        &    5 &    2 &      ---       &   ---    &   ---    & 1.14e+03 & 10 \\
BT8        &    4 &    2 &  +1.000243e+00 & 1.98e-03 & 3.41e-04 & 1.46e+01 & 10 \\
BT9        &    2 &    2 &  -1.000035e+00 & 1.17e-02 & 2.21e-04 & 1.98e+02 & 10 \\
BT10       &    5 &    3 &  -1.000038e+00 & 0.00e+00 & 5.91e-05 & 5.00e+00 & 10 \\
BT11       &    5 &    3 &  +8.258014e-01 & 1.25e-01 & 5.58e-04 &    ---   &  5 \\
BT12       &    3 &    2 &  +6.189020e+00 & 2.82e-02 & 3.48e-04 & 5.23e+02 & 10 \\
BYRDSPHR   &   10 &    5 &  -4.664093e+00 & 1.11e-01 & 1.19e-04 & 8.67e+02 & 10 \\
DIXCHLNG   &  110 &   55 &  +2.472434e+03 & 6.46e+01 & 5.25e-04 &    ---   &  5 \\
EIGENA2    &  110 &   55 &  +9.973427e-07 & 1.86e-03 & 3.57e-04 & 4.15e+02 & 10 \\
EIGENACO   &  110 &   55 &  +2.948147e-07 & 1.08e-03 & 9.41e-16 & 2.81e+02 & 10 \\
EIGENB2    &  110 &   55 &  +3.675137e-06 & 1.07e-03 & 5.37e-04 & 8.56e+02 & 10 \\
EIGENBCO   &   75 &   25 &  +2.125972e-06 & 7.45e-04 & 5.15e-04 & 1.94e+03 & 10 \\
ELEC       &   10 &    8 &  +2.438199e+02 & 6.84e-01 & 5.62e-04 &    ---   &  2 \\
GENHS28    &    7 &    2 &  +9.280133e-01 & 5.07e-02 & 6.75e-06 &    ---   &  1 \\
HS100LNP   &    2 &    1 &  +6.806891e+02 & 1.84e+00 & 8.94e-04 &    ---   &  0 \\
HS6        &    2 &    1 &  +1.921030e-03 & 1.60e-02 & 3.73e-04 & 2.96e+01 & 10 \\
HS7        &    2 &    2 &  -1.732161e+00 & 8.73e-04 & 3.80e-04 & 9.28e+01 & 10 \\
HS8        &    2 &    1 &  -1.000000e+00 & 0.00e+00 & 1.32e-08 & 4.00e+00 & 10 \\
HS9        &    3 &    1 &  -4.968118e-01 & 1.36e-02 & 1.10e-14 & 7.22e+01 & 10 \\
HS26       &    3 &    1 &  +3.827215e-05 & 2.15e-03 & 4.57e-04 & 1.67e+02 & 10 \\
HS27       &    3 &    1 &  +4.007659e-02 & 2.90e-03 & 2.31e-04 & 9.16e+01 & 10 \\
HS28       &    4 &    2 &  +2.028082e-05 & 2.71e-03 & 8.44e-16 & 7.52e+01 & 10 \\
HS39       &    4 &    3 &  -1.000374e+00 & 5.82e-03 & 4.58e-04 & 1.94e+02 & 10 \\
HS40       &    4 &    2 &  -1.981579e-01 & 4.81e-02 & 1.35e-04 & 3.81e+02 & 10 \\
HS42       &    5 &    2 &  +1.391016e+01 & 5.52e-01 & 2.13e-04 & 4.87e+03 & 10 \\
HS46       &    5 &    3 &  +4.013015e-05 & 1.67e-03 & 5.27e-04 & 1.56e+02 & 10 \\
HS47       &    5 &    2 &  +5.386687e-05 & 2.46e-03 & 5.13e-04 & 4.49e+02 & 10 \\
HS48       &    5 &    3 &  +8.988756e-06 & 4.31e-03 & 1.39e-15 & 8.36e+01 & 10 \\
HS50       &    5 &    3 &  +8.049833e-05 & 1.25e-02 & 1.08e-13 & 5.38e+03 & 10 \\
HS51       &    5 &    3 &  +2.340881e-07 & 9.39e-04 & 1.98e-15 & 1.68e+01 & 10 \\
HS52       &    3 &    2 &  +5.328581e+00 & 1.57e-01 & 6.34e-05 &    ---   &  0 \\
HS61       &    5 &    2 &      ---       &   ---    &   ---    & 1.00e+00 & 10 \\
HS77       &    5 &    3 &  +2.416904e-01 & 3.22e-02 & 4.66e-04 &    ---   &  0 \\
HS78       &    5 &    3 &  -2.914874e+00 & 1.41e-01 & 2.74e-04 &    ---   &  0 \\
HS79       &   20 &   18 &  +7.880086e-02 & 4.41e-03 & 5.19e-04 & 9.58e+02 & 10 \\
LUKVLE1    &   20 &   13 &  +3.740180e+00 & 1.06e-02 & 2.27e-04 & 1.11e+02 & 10 \\
LUKVLE2    &   20 &    2 &  +4.371965e+02 & 8.04e+00 & 5.47e-04 &    ---   &  6 \\
LUKVLE3    &   20 &    9 &  +2.758515e+01 & 1.13e+00 & 4.96e-04 &    ---   &  5 \\
LUKVLE4    &   21 &   10 &  +1.859096e+02 & 2.56e+02 & 9.65e-05 & 1.00e+05 & 10 \\
LUKVLE6    &   20 &    4 &  +1.119752e+03 & 4.22e+00 & 5.48e-04 & 1.00e+05 & 10 \\
LUKVLE7    &   20 &   18 &  -5.347145e+00 & 5.86e-01 & 4.82e-04 &    ---   &  0 \\
LUKVLE8    &   20 &    6 &      ---       &   ---    &   ---    & 1.90e+02 & 10 \\
LUKVLE9    &   20 &   18 &  +7.100923e+00 & 1.90e+01 & 3.53e-04 & 1.00e+05 & 10 \\
LUKVLE10   &   18 &   10 &  +6.459921e+00 & 2.75e-02 & 4.97e-04 & 9.69e+04 & 10 \\
LUKVLE11   &   17 &   12 &  +1.000256e-04 & 2.75e-03 & 5.29e-04 & 3.72e+03 & 10 \\
LUKVLE12   &   18 &   10 &  +2.286872e+02 & 4.44e+01 & 2.23e+00 & 1.00e+05 & 10 \\
LUKVLE13   &   18 &   10 &  +5.501757e+01 & 3.78e-01 & 3.38e-04 & 1.00e+05 & 10 \\
LUKVLE14   &   17 &   12 &  +4.180692e+04 & 4.10e+00 & 3.72e-04 & 1.00e+05 & 10 \\
LUKVLE15   &   17 &   12 &  +7.963970e+01 & 1.27e+01 & 1.81e-04 &    ---   &  0 \\
LUKVLE16   &   17 &   12 &  +8.040898e+01 & 2.71e+01 & 3.79e+00 &    ---   &  0 \\
LUKVLE17   &   17 &   12 &  +1.190180e+02 & 5.48e+01 & 8.37e+00 &    ---   &  0 \\
LUKVLE18   &   20 &    9 &  +2.627594e+01 & 1.71e+01 & 8.37e+00 &    ---   &  0 \\
LUKVLI4    &    2 &    1 &  +1.641033e+02 & 2.16e+02 & 2.26e-04 & 1.00e+05 & 10 \\
MARATOS    &    5 &    3 &  -1.000135e+00 & 4.39e-03 & 3.68e-04 & 3.51e+01 & 10 \\
MWRIGHT    &    9 &    3 &  +2.498269e+01 & 4.85e-01 & 4.89e-04 &    ---   &  4 \\
ORTHRDM2   &    9 &    3 &  +1.938152e-07 & 7.12e-04 & 2.85e-04 & 2.68e+01 & 10 \\
ORTHRDS2   &  133 &   64 &  +2.075502e-07 & 7.50e-04 & 2.86e-04 & 2.67e+01 & 10 \\
ORTHREGA   &   27 &    6 &  +3.503053e+02 & 2.04e-01 & 6.55e-04 &    ---   &  0 \\
ORTHREGB   &   15 &    5 &  +1.903767e-07 & 7.21e-04 & 2.66e-04 & 1.51e+01 & 10 \\
ORTHREGC   &   43 &   20 &  +7.509051e-07 & 8.31e-04 & 5.59e-04 & 3.67e+01 & 10 \\
ORTHREGD   &   43 &   20 &      ---       &   ---    &   ---    & 1.65e+02 & 10 \\
ORTHRGDM   &   43 &   20 &      ---       &   ---    &   ---    & 1.84e+03 & 10 \\
ORTHRGDS   &    2 &    1 &  +6.212538e+00 & 9.40e-02 & 7.00e-04 &    ---   &  0 \\
S316m322   &   11 &    9 &      ---       &   ---    &   ---    & 0.00e+00 & 10 \\
SPINOP     &   11 &    9 &  +3.658931e-01 & 1.59e-02 & 4.10e-04 & 1.53e+03 & 10 \\
\hline
\end{tabular}
\end{center}
}
\end{table}

\end{document}